\documentclass[a4paper,11pt,english]{article}

\usepackage{amssymb}
\usepackage{mathtools}    
\usepackage[all]{xy}
\usepackage{amsthm}
\usepackage{currvita}
\usepackage{mathdots}
\usepackage{amsmath}
\usepackage{verbatim}
\usepackage{MnSymbol}
\usepackage{makeidx}
\usepackage{mathrsfs}
\usepackage[ansinew]{inputenc}
\usepackage{microtype}
\SelectTips{cm}{12}
\makeindex
\newtheorem{thm}{Theorem}[section]
\newtheorem{lemma}[thm]{Lemma}
\newtheorem{defi}[thm]{Definition}
\newtheorem{defilemma}[thm]{Definition/Lemma}
\newtheorem{cor}[thm]{Corollary}
\newtheorem{rem}[thm]{Remark}
\newtheorem{prop}[thm]{Proposition}
\newtheorem{exa}[thm]{Example}

\renewcommand{\ker}{\operatorname{ker}}

\newcommand{\id}{\operatorname{id}}
\newcommand{\im}{\operatorname{im}}
\newcommand{\lin}{\operatorname{lin}}
\newcommand{\rank}{\operatorname{rank}}

\newcommand{\sgn}{\operatorname{sgn}}
\newcommand{\vp}{{\varphi}}

\newcommand{\diago}{\operatorname{diag}}

\newcommand{\Mult}{\operatorname{Mult}}

\newcommand{\cR}{{\mathcal R}}

\newcommand{\cL}{{\mathcal L}}

\newcommand{\N}{{\mathbb N}}
\newcommand{\K}{{\mathbb K}}
\renewcommand{\O}{{\mathbb O}}
\newcommand{\lr}[1] {{\left(#1\right)}}
\newcommand{\fr}[1] {{\frac{1}{#1}}}

\newcommand{\set}[2]{\left\{#1\left|\; #2\right.\right\}}

\newcommand{\bt}{\begin{thm}}
\newcommand{\et}{\end{thm}}
\newcommand{\bp}{\begin{proof}}
\newcommand{\ep}{\end{proof}}
\newcommand{\bl}{\begin{lemma}}
\newcommand{\el}{\end{lemma}}
\newcommand{\bc}{\begin{cor}}
\newcommand{\ec}{\end{cor}}
\newcommand{\bd}{\begin{description}}
\newcommand{\ed}{\end{description}}

\newcommand{\fA}{{\mathfrak A}}
\newcommand{\fB}{{\mathfrak B}}

\newcommand{\fG}{{\mathfrak G}}

\newcommand{\fg}{\mathfrak g}

\newcommand{\M}{\mathbb M}
\newcommand{\Z}{\mathbb Z}
\newcommand{\Li}{\mathbb L}

\newcommand{\cG}{\mathcal G}

\newcommand{\sub}{\subseteq}

\begin{document}

\thispagestyle{empty}
\begin{center}
\huge
\textsf{A K-theoretic approach to the classification of symmetric spaces}\\
\vspace*{1cm} \large
Dennis Bohle, Wend Werner\\
Fachbereich Mathematik und Informatik\\
Westf\"alische Wilhelms-Universit\"at\\
Einsteinstra\ss e 62\\
48149 M\"unster \\
\vspace*{0.5cm}
 dennis.bohle@math.uni-muenster.de,\\ wwerner@math.uni-muenster.de\\
\vspace*{0.3cm}
\end{center}
\begin{abstract}

We show that the classification of the
symmetric spaces can be achieved by K-theoretical methods. We focus on
Hermitian symmetric spaces of non-compact type, and define K-theory
for JB*-triples along the lines of C*-theory. K-groups have to be provided
with further invariants in order to classify.
Among these are the cycles obtained from so called \emph{grids}, intimately
connected to the root systems of an underlying Lie-algebra and thus reminiscent
of the classical classification scheme.

\footnotetext{2000
\emph{Mathematics Subject Classification}: 17C65, 46L70}

\footnotetext{Key words and phrases: $JB^*$-triple system, K-theory,
Hermitian symmetric space, grid, Cartan factor, TRO, ternary ring of operators,
 universal enveloping TRO} \footnotetext{The first
author was supported by the Graduiertenkolleg f\"ur analytische
Topologie und Metageometrie}

\end{abstract}

In the following we are going to prove that Cartan's classification of the
symmetric spaces has a K-theoretic background. We will do so within the
realm
of Hermitian symmetric spaces of non-compact type; the more general
Riemannian symmetric spaces are either within reach by passing to duals
(in the compact case), or, respectively, should be classifiable by real
K-theory,
still to be developed.

By its Harish-Chandra realization and the Koecher theory, Hermitian
symmetric
spaces of non-compact type allow a realization as the open unit ball of
a JB*-triple
system. The latter generalize C*-algebras and sport a triple product which
resembles in its basic properties their binary C*-counterparts.

The idea behind our classification scheme is to define for the category
of JB*-triple
systems a K-functor along the lines of C*-algebra theory and then
provide the emerging
groups with further invariants, sharp enough to distinguish (at least)
the finite
dimensional JB*-triple systems.

One of the major impediments along this way is the fact that JB*-triple
systems do
not, in general, behave well under the formation of tensor products.
Those which do,
coincide almost exactly with the Hilbert-C*-modules, and for this
category we define K-theory
in the second section (after we have collected some preliminary
information in the first).
In order to emphasize their relationship to the more general JB*-triple
systems (and, simultaneously, have the right morphisms at hand) we
treat these modules in the guise of \emph{ternary rings of operators
(TROs)}.

The functor from section 2 is concatenated with the one from
\cite{BoWe1, BunceFeelyTimoneyI},
which assigns a universal TRO to any JB*-triple system. Both functors
are shown to behave
so nicely, that they yield K-theory with all expected properties,
except for stability, which,
as noted, is notoriously hard to come by right from the start. An
interesting point is the
fact that a C*-algebra, equipped with its natural JB*-structure, has in
general a K-theory
within the context of JB*-triples which differs from its C*-version.

As in the case of C*-algebras, K-groups themselves do not suffice to
distinguish between
different objects of interest. We therefore enhance the K-groups for the
JB*-triple systems
by a twofold order structure ("scales") as well as by a subset of
distinguished cycles.
The latter come from so called  \emph{grids}, which form a set of
special generators and
relations. That these objects in combination form an invariant
(which we will call the
K-grid invariant) is shown in section 4.
We calculate the K-grid invariant for the Cartan factors in section 5.
As it turns out, it discriminates all but the
exceptional factors (for which the K-groups vanish).

The results of this paper are taken from the first named author's PhD-thesis
\cite{Bohle}.

\section{Preliminaries}\label{section Preliminaries}
Let us begin by recalling some of the more general concepts that
will be used in the following.

A generalization of C* -algebras are the JB*-triple systems, or JB*-triples for
short. Instead of the binary product and the involution of the former, they come
equipped with  a ternary product $(a,b,c)\mapsto\{a,b,c\}$ that is symmetric
complex bilinear in the outer factors $a, c$ and conjugate linear in the middle
factor $b$. Every C*-algebra, for instance, is a JB*-triple with respect to the
Jordan triple product $\{a,b,c\} = (ab^*c+cb^*a)/2$.
For the exact definition of JB*-triples
(as well as a coherent presentation of the theory) we refer the reader to
\cite{Upmeier-SymmetricBanachmanifoldandJordaCalgebras}.

In finite dimension, any JB*-triple can be written as the sum of so called
\emph{Cartan factors}, the prime elements of the theory in finite dimensions
\cite{FriedmanRusso-TheGelfandNaimarktheoremforJBtriples}.
There are four series,
\begin{itemize}
\item
the \emph{rectangular factors} $\M_{n,m}$ of $n\times m$-matrices,
\item
\emph{hermitian factors}, consisting of the $n\times n$ symmetric matrices,
\item
the \emph{symplectic factors}, which are formed by the $n\times n$
skew-symmetric matrices,
\item
the spin factors, which are linearly generated by the identity on an $n$-dimensional
complex Hilbert space together with $n$ self-adjoint matrices $s_1,\ldots,s_n$
fulfilling the canonical anti-commutator relation $\frac{1}{2}(s_i s_j+s_j s_i)=\delta_{i,j}$,
\end{itemize}
as well two exceptional ones,
\begin{itemize}
\item
the space of all 1 x 2 matrices over $\O$, the complex Cayley numbers
(with triple product  $\{x,y,z\} = 1/2[x(y^*z)+z (y^*x)]$), and
\item
the  space of all 3 x 3 hermitian matrices over $\O$. The triple
product here is defined via the Jordan product
$\{x,y,z\}=(x\circ y)\circ z+(z\circ y)\circ x-(x\circ z)\circ y$,
where $x\circ y = 1/2 (xy+yz)$.
\end{itemize}

The latter two distinguish themselves by not being isomorphic to a JB*-subtriple
of $B(H,K)$, the space of bounded operators between Hilbert spaces $H$ and $K$.

Further examples are given by the class of (full) Hilbert-C*-modules,
for which the triple product in terms of the left and right inner products $\langle\cdot,
\cdot\rangle_r$ and $\langle\cdot,\cdot\rangle_l$ are given by
$$
\{a,b,c\}=\frac{1}{2}\left(\langle a,b\rangle_lc+a\langle b,c\rangle_r\right).
$$
Actually, this class of Hilbert-C*-modules coincides with the so called
\emph{ternary rings of operators}, TROs for short. These are the closed
subspaces $T$ of  B(H,K) such that the ternary product
$xy^*z$ lies in $T$ for all $x,y,z \in T$. If, given a TRO, we identify
$$
\cL(T)=TT^*=\overline{\lin}\set{xy^*}{x,y\in T}
\quad\text{and}\quad
\cR(T)=T^*T=\overline{\lin}\set{x^*y}{x,y\in T}
$$
with a subalgebra of bounded linear operators on $T$, then both these spaces are
C*-algebras (called the \emph{left} and, respectively, \emph{right
C*-algebra of $T$}), and $T$ becomes a (two-sided) Hilbert-C*-module over
$\cR(T)$ and $\cL(T)$.

Here are three important features to be kept in mind:
\begin{itemize}
\item
Other than for an arbitrary JB*-triple, the space
of matrices $M_n(T)$ with entries from T again is a TRO, in the canonical way.
This property actually characterizes TROs among the JB*-triples
\cite{Hamana-TripleenvelopesandSilovboundariesofoperatorspaces}
\item
The morphisms for these objects, TRO-morphisms, are different from those that
are generated by the underlying Hilbert-C*-module.
\item
For each TRO $T$ there exists a canonical C*-algebra, the \emph{linking algebra $\Li(T)$},
formally defined through
$$
\Li(T)=\begin{pmatrix}\cR(T)& T\\
                        T^* &\cL(T)\end{pmatrix},
$$
with $T^*$ the TRO conjugate to $T$.
\end{itemize}
The expressions $\cL,\cR$, and $\Li$ actually become functorial by letting
\begin{gather*}
\cL(\Phi)(\sum x_iy_i^*)=\sum \Phi(x_i)(y_i)^*\qquad
\cR(\Phi)(\sum x_i^*y_i)=\sum \Phi(x_i)^*(y_i)\\
\text{and}\qquad
\Li(\Phi)=\begin{pmatrix}\cR(\Phi)& \Phi\\
                        \Phi^* &\cL(\Phi)\end{pmatrix}
\end{gather*}
for any TRO-morphism $\Phi$. A standard reference for these topics is \cite{BlecherLeMerdy-Operatoralgebrasandtheirmodules}

Finally, we give a brief account of an important link between JB*-theory and the classical
theory of symmetric
spaces, which essentially was developed by E.\ Neher \cite{Neher-Jordantriplesystemsbythegridapproach,
Neher-Systemesderacines3gradues,
Neher-3gradedrootsystemsandgridsinJordantriplesystems}.
Each of the Cartan factors can be defined by a special set of generators and relations,
the grids. These generators consist of tripotent elements $u$, $\{u,u,u\}=u$ and are related to root systems, in the following way.
For each JB*-triple let $Z_{1}=Z_{-1}=Z$ and $\partial Z=\set{z\mapsto\{a,b,z\}}{a,b\in Z}$. Then
$$
\fg(Z)=Z_{-1}\oplus\partial Z\oplus Z_1
$$
can be equipped with a Lie-bracket and is called the \emph{Tits-Kantor-Koecher Lie algebra (or TKK-algebra) of $Z$}. This algebra is 3-graded in an obvious fashion. There exists
furthermore a theory of 3-graded root systems, where a root system $\cR$ is decomposed into a disjoint union $\cR=\cR_{-1}\cup\cR_0\cup\cR_1$. Then, for each of the Cartan factors, there
is a 3-graded root system, and for each such there is associated a grid with its $\cR_1$ part
\cite{Neher-3gradedrootsystemsandgridsinJordantriplesystems}.

\section{K-theory for TROs}

\subsection{The Linking Functor}

In this section we show that $\cL,\cR$ and $\Li$ are covariant functors from the category of TROs with TRO-homomorphisms to the category of $C^*$-algebras with $*$-homomorphisms. We determine explicitly the properties of these functors. All three functors are exact, homotopy invariant, stable and continuous and therefore especially additive and split exact which makes them all excellent candidates to define a functor $K_0$ from the category of TROs to the category of Abelian groups. We choose the functor $\cL$, but this does not affect the general theory since all resulting $K_0$-groups are isomorphic. For the sake of brevity let $\kappa\in\{\cL,\cR,\Li\}$.

It is easy to see that the mappings $\cL$, $\cR$ and $\Li$ are covariant functors from the category of TROs to the category of $C^*$-algebras.
\begin{defi}
Let $T$ and $U$ be TROs. Two TRO-homomorphisms
$\varphi,\psi:T\rightarrow U$ are called \textbf{homotopic}\index{homotopic TROs} (denoted by $\varphi\sim_h
\psi$), when there exists a path of TRO-homomorphisms
$\gamma_t:T\rightarrow U$, $t\in[0,1]$, such that $t\mapsto
\gamma_t(x)$ is a continuous map from $[0,1]$ to $U$ for all $x\in T$, satisfying  $\gamma_0=\varphi$ and
$\gamma_1=\psi$. The TROs $T$ and $U$ are called \textbf{homotopy equivalent}\index{homotopy equivalent TROs} if there are TRO-homomorphisms $\varphi:T\rightarrow U$ and
$\psi: U\rightarrow T$ with $\varphi\circ\psi\sim_h \id_U$
and $\psi \circ\varphi\sim_h \id_T$.
\end{defi}

The proof of the next proposition is obvious and thus omitted.

\begin{prop}\label{homotopic tro homomorphisms}
Let $T$ and $U$ be TROs.
\begin{description}
\item[(a)] If $\varphi,\psi:T\rightarrow U$ are homotopic TRO-homomorphisms, then $\kappa(\varphi)$ and
$\kappa(\psi)$ are homotopic $*$-homomorphisms.
\item[(b)] If $T$ and $U$ are homotopic TROs, then $\kappa(T)$
and $\kappa(U)$ are homotopic $C^*$-algebras.
\end{description}
\end{prop}

\begin{cor}
Two $C^*$-algebras $\mathfrak{A}$ and $\mathfrak{B}$ are homotopic as  $C^*$-algebras if and only if they are homotopic as TROs.
\begin{proof}
If $\fA$ and $\fB$ are homotopic as $C^*$-algebras they are homotopic as TROs, since every $*$-homomorphism is a TRO-homomorphism.

Using Proposition \ref{homotopic tro homomorphisms} (b) we get that $\cL(\fA)=\fA$ and $\cL(\fB)=\fB$ are homotopic as $C^*$-algebras.
\end{proof}
\end{cor}

That our functors are exact and therefore in particular split exact and additive follows from the Rieffel quotient equivalence (c.f.\ \cite{BlecherLeMerdy-Operatoralgebrasandtheirmodules}, 8.2.25).

To show that inductive limits exist in the category of TROs and that the linking algebra of an inductive limit of TROs equals the inductive limit of the corresponding linking algebras, we first have to consider the following special case, the proof of which is standard.

\begin{lemma}
Let $T$ be a TRO and $(T_n)$ an increasing sequence of subTROs of $T$. Denote by $\vp_n:T_n\to T_{n+1}$ the inclusion mapping and put $$T_\infty:=\overline{\bigcup_{n=1}^{\infty}T_n}.$$ Then $(T_\infty,(\iota_n))$ is the inductive limit of $((T_n),(\vp_n))$, where $\iota_n$ denotes the inclusion mapping of $T_n$ into $T_\infty$ for all $n\in\N$.
\end{lemma}

Next we show the existence of inductive limits in the category of TROs.

\begin{prop}\label{inductive lim fuer TROs}
Every inductive system $((T_n),(\varphi_n))$ in the category of TROs has an inductive limit $(T_\infty,(\mu_n))$.
\begin{proof}
We first note that
 $((\mathcal{L}(T_n)),(\mathcal{L}(\varphi_n)))$, $((\mathcal{R}(T_n)),(\mathcal{R}(\varphi_n)))$ and \linebreak $((\Li(T_n)),$ $(\Li(\varphi_n)))$ are inductive systems of
$C^*$-algebras, whose inductive limits exist by $C^*$-therory. Let $(\Li_\infty,(\lambda_n))$ be the inductive limit of the sequence of linking algebras then by \cite{RordamLarsenLaustsen-AnintroductiontoKtheoryforCalgebras}, 6.2.4 (i),
\[\Li_\infty=\overline{\bigcup_{n=1}^\infty\lambda_n\left( \Li(T_n)\right)}.\]
Let $\iota_n:T_n\to\Li(T_n)$ be the canonical corner embeddings and\[
T_\infty:=\overline{\bigcup_{n=1}^\infty\lambda_n\left( \iota_n(T_n)\right)}\sub\Li_\infty.
\]
Define $\mu_n:=\lambda_n\circ\iota_n:T_n\to T_\infty$, then
\begin{align*}
\mu_{n+1}\circ\vp_n(x)&=\lambda_{n+1}\circ\iota_{n+1}\circ\vp_n(x)\\
&=\lambda_{n+1}\circ\Li(\varphi_n)\begin{pmatrix}
0&x\\
0&0
\end{pmatrix}\\
&=\lambda_n\circ\iota_n(x)\\
&=\mu_n(x)
\end{align*}
for all $x\in T_n$, $n\in\N$. All $\mu_n$ are TRO-homomorphisms and $T_\infty$ is a subTRO of $\Li_{\infty}$.

To prove the uniqueness of the inductive limit let $(U,(\beta_n))$ be another system satisfying $\beta_{n+1}\circ\vp_n=\beta_n$, where $\beta_n:T_n\to U$ is a TRO-homomorphism for all $n\in \N$. Since $(\Li_\infty,(\lambda_n))$ is the inductive limit of $(\Li(T_n),(\Li(\varphi_n)))$ and $\Li(\beta_{n+1})\circ\Li(\vp_n)=\Li(\beta_{n+1}\circ\vp_n)=\Li(\beta_n)$, there exists one and only one $*$-homomorphism $\lambda$ making the diagram
\begin{equation*}
\text{\xymatrix{
 &\Li(T_n)\ar[dl]_{\lambda_n}\ar[dr]^{\Li(\beta_{n})}
\\
\Li_\infty
\ar[rr]^-{\lambda}
&&
\Li(U)
 }}\end{equation*} commutative. The restriction of $\lambda$ to $T_\infty$ gives the desired TRO-homo\-morphism from $T_\infty$ to $U$.
\end{proof}
\end{prop}

To prove that our functors are continuous we need the following lemma.

\begin{lemma}\label{hilfefurindlimiten}
Let $((T_n),(\vp_n))$ be an inductive sequence of TROs with inductive limit $(T_\infty,(\mu_n))$, then

$$\ker(\mu_n)=\{x\in T_n: \lim_{m\to\infty}||\vp_{m,n}(x)||=0\}$$
for all $n\in\N$.

\begin{proof} Let $(\Li_\infty,(\lambda_n))$ be the inductive limit of the sequence of linking algebras $(\Li(T_n),(\vp_n))$.
We know from the proof of Proposition \ref{inductive lim fuer TROs} that $\Li_\infty=\overline{\bigcup_{n=1}^{\infty}\lambda_n(\Li(T_n))}$,
$T_\infty=\overline{\bigcup_{n=1}^{\infty}\lambda_n(\iota_{T_n}(T_n))}\sub\cL_\infty$ and $\mu_n=\lambda_n\circ\iota_{T_n}$ for all $n\in\N$. Since $(\Li_\infty,(\lambda_n))$ is a $C^*$-direct limit we know from \cite{RordamLarsenLaustsen-AnintroductiontoKtheoryforCalgebras}, Proposition 6.2.4 that $$||\lambda_n(x)||=\lim_{m\to\infty}||\Li(\vp_{m,n})(x)||$$ for all $n\in \N$. Especially for all $y\in T_n$ we get
\begin{align*}
||\mu_n(y)||&=||\lambda_n(\iota_{T_n}(y))||\\
&=\lim_{m\to\infty}||\Li(\vp_{m,n})(\iota_{T_n}(y))||\\
&=\lim_{m\to\infty}\left|\left|\begin{pmatrix}0&\vp_{m,n}(y)\\
0&0\end{pmatrix}\right|\right|\\
&=\lim_{m\to\infty}\left|\left|\vp_{m,n}(y)\right|\right|.
\end{align*}
\end{proof}
\end{lemma}

We can now show that the functors $\cL$, $\cR$ and $\Li$ are continuous.

\begin{prop}
If $((T_n),(\vp_n))$ is an inductive sequence of TROs with inductive limit $(T_\infty,(\mu_n))$, then \[\lim_{n\to\infty}\kappa(T_n)=\kappa(T_\infty).\]
\begin{proof}
First recall from the proof of Proposition \ref{inductive lim fuer TROs} that the inductive limit of the inductive sequence $((\Li(T_n),(\Li(\vp_n)))$ is $\Li_\infty=\overline{\bigcup_{n=1}^\infty\lambda_n(\Li(T_n))}$ and that $T_\infty=\overline{\bigcup_{n=1}^\infty\mu_n(T_n)}$, with $\mu_n:=\lambda_n\circ\iota_{T_n}:T_n\to T_\infty$. If we put $$\cL_\infty:=\overline{\bigcup_{n=1}^\infty\lambda_n(\iota_{\cL(T_n)}(\cL(T_n))}\sub\Li_\infty$$ we see that this is the inductive limit of the sequence $((\cL(T_n)),\cL(\vp_n))$ and
\begin{align*}
\cL(T_\infty)
&=\overline{\bigcup_{n=1}^\infty\left(\lambda_n(\iota_{T_n}(T_n))\right)\left(\lambda_n(\iota_{T_n}(T_n))\right)^*}
=\cL_\infty,
\end{align*}
which shows that $\lim_{n\to\infty}\cL(T_n)=\cL(T_\infty)$. The homomorphisms $\eta_n:\cL(T_n)\to\cL(T_\infty)$ are given by $\eta_n:=\cL(\mu_n)=\cL(\lambda_n\circ\iota_{T_n})=\lambda_n\circ\iota_{\cL(T_n)}$ for all $n\in\N$. A similar proof shows that $\left(\cR(T_\infty),(\cR(\mu_n))\right)$ is the inductive limit of $\left((\cR(T_n)),(\cR(\vp_n))\right)$.

Finally we show that $\Li_\infty=\lim_{n\to\infty}\Li(T_n)=\Li(T_\infty)$ and that $\lambda_n=\Li(\mu_n)$ for all $n\in\N$. Since $(T_\infty,(\mu_n))$ is the inductive limit of $((T_n),(\vp_n))$ we notice that
$\mu_{n+1}\circ\vp_n=\mu_n$
and thus $$\Li(\mu_{n+1})\circ\Li(\vp_n)=\Li(\mu_n)$$ for all $n\in\N$.

We get by the universal property of the inductive limit a unique $*$-homomorphism $\lambda:\Li_\infty\to\Li(T_\infty)$ and for all $n\in\N$ a commutative diagram
\[
\text{\xymatrix{
 &\Li(T_n)\ar[dr]^{\Li(\mu_n)}\ar[dl]_{\lambda_{n}}
\\
\Li_{\infty}
\ar[rr]^-{\lambda}
&&
\Li(T_\infty)
 }}
\]

It is well known for inductive limits of $C^*$-algebras (cf.\ \cite{RordamLarsenLaustsen-AnintroductiontoKtheoryforCalgebras}, Proposition 6.2.4), that
\[\lambda \text{ is surjective }\Longleftrightarrow \Li(T_\infty)=\overline{\bigcup_{n=1}^\infty\Li(\mu_n)\Li(T_n)}, \] which is the case, and that
\[\lambda \text{ is injective }\Longleftrightarrow \ker(\Li(\mu_n))\sub\ker(\lambda_n)\]
for all $n\in\N$. It follows from \cite{RordamLarsenLaustsen-AnintroductiontoKtheoryforCalgebras}, Proposition 6.2.4 and the above that \begin{align*}
\ker(\lambda_n)&=\{x\in \Li(T_n): \lim_{m\to\infty}||\Li(\vp_{m,n})(x)||=0\},\\
\ker(\cR(\vp_n))&=\{x\in \cR(T_n): \lim_{m\to\infty}||\cR(\vp_{m,n})(x)||=0\}
 \end{align*}
 and
 \begin{align*}
\ker(\cL(\vp_n))=\{x\in \cL(T_n): \lim_{m\to\infty}||\cL(\vp_{m,n})(x)||=0\}
 \end{align*}
 for all $n\in\N$.
 Since $(T_\infty,(\mu_n))$ is the inductive limit of $((T_n),(\vp_n))$ we get with Lemma \ref{hilfefurindlimiten} that
$$\ker(\mu_n)=\{y\in T_n: \lim_{m\to\infty}||\vp_{m,n}(y)||=0\},$$
for all $n\in\N$.

Let $n\in\N$ and $x=\begin{pmatrix}x_1&x_2\\x_3&x_4\end{pmatrix}\in\ker(\Li(\mu_n))\sub\Li(T_n)$,
then $$\lim_{m\to\infty}||\vp_{m,n}(x_2)||=0$$ and by $C^*$-theory $$\lim_{m\to\infty}||\cL(\vp_{m,n})(x_4)||=0=\lim_{m\to\infty}||\cR(\vp_{m,n})(x_1)||.$$ Thus we get
 \begin{align*}\lim_{m\to\infty}\left|\left|\Li(\vp_{m,n})(x)\right|\right|&= \lim_{m\to\infty}\left|\left|\begin{pmatrix}\cR(\vp_{m,n})(x_1)&\vp_{m,n}(x_2)\displaybreak[0]\\
\vp_{m,n}^*(x_3)&\cL(\vp_{m,n})(x_4)
\end{pmatrix}\right|\right|\displaybreak[0]\\
&\leq\lim_{m\to\infty}(||\cR(\vp_{m,n})(x_1)||+||\vp_{m,n}(x_2)||\displaybreak[0]\\
&{}+||\vp_{m,n}^*(x_3)||+||\cL(\vp_{m,n})(x_4)||)\displaybreak[0]\\
&=0.
\end{align*}
Therefore we can conclude that $x\in\ker{\lambda_n}$ and that $\lambda$ is an isomorphism.
\end{proof}
\end{prop}

\begin{prop}
If $T$ is a TRO, then
\begin{align*}\kappa(M_n(T))&\simeq M_n(\kappa(T))
\end{align*}
are isomorphic as $C^*$-algebras.
\begin{proof}
 This follows easily from fact that the permutation of rows and columns yields complete isometries.
\end{proof}
\end{prop}

\begin{defilemma}
Let $T$ be TRO. We define the \textbf{(ternary)} \textbf{K-groups}\index{ternary K-groups} of $T$ to be
\begin{align*}K_i^{\text{\tiny TRO}}(T):=K_i^{\text{\tiny C*}}(\mathcal{L}(T))\text{ for }i\in \N_0.
\end{align*}

Every TRO-homomorphism $\varphi: T\to U$ induces a homomorphism of Abelian groups\[
K_i^{\text{\tiny TRO}}(\varphi):K_i^{\text{\tiny TRO}}(T)\to K_i^{\text{\tiny TRO}}(U),
\] when we put $K^{\text{\tiny TRO}}_{i}(\vp):=K_i^{\text{\tiny C*}}(\mathcal{L}(\vp))$ for $i\in\N_0$.
\end{defilemma}

The ternary $K$-groups of a $C^*$-algebra $\fA$ are exactly the $K$-groups in the $C^*$-algebra sense, since $\fA\fA^*=\fA$. Therefore it poses no danger of confusion to write $K_i(\fA)$ for the ternary $K_i$-group of $\fA$, for $i\in\N_0$. If $\vp: \fA\to\fB$ is a $*$-homomorphism we know that $\mathcal{L}(\vp)=\vp$ and thus $K_i^{\text{\tiny TRO}}(\vp)=K_i^{\text{\tiny C*}}(\mathcal{L}(\vp))=K_i^{\text{\tiny C*}}(\vp)$. We therefore drop the superscripts for $K_i$.

\begin{prop}\label{K_0 funktor eigenschaften}
Let $i\in\N_0$. The functors $K_i$ are half-exact, split-exact, homotopy invariant, continuous and Bott periodic.

\begin{proof}
This is just basic $K$-theory combined with our functorial results for $\cL$ from this section.
\end{proof}
\end{prop}

The next proposition shows that the choice we made by preferring the functor $\cL$ over the functors $\cR$ and $\Li$ does not affect the theory. To prove this theorem we make use of the theory of stably isomorphic $C^*$-algebras developed by Brown. Since this theory is intended for separable $C^*$-algebras, we restrict our attention to separable TROs.

\begin{prop}\label{moritaisosatz}
Let $T$ be a separable TRO, then \[K_0(T)=K_0(\mathcal{L}(T))\simeq K_0(\mathcal{R}(T))\simeq K_0(\Li(T)).\]
\begin{proof} It follows from \cite{Brown-StableisomorphismofhereditarysubalgebrasofCalgebras}, Lemma 2.5 and Lemma 2.6 that, if $\fA$ is a $C^*$-algebra containing a strictly positive element and $\mathfrak{B}$ is a full corner of $\fA$ (i.e.\ there exists a projection $p$ in the multiplier algebra $\Mult(\fA)$ of $\fA$ such that $\mathfrak{B}=p\fA p$ is not contained in any proper closed two-sided ideal of $\fA$), then there exists a tripotent $v$ in the multiplier algebra of $\fA\otimes \K$ such that $v^*v=\id$ and $vv^*=p\otimes\id$. Moreover, $\fA$ and $\fB$ are stably isomorphic and the isomorphism is induced by the partial isometry $v$, here $p$ is the projection with $\mathfrak{B}=p\fA p$ and $\mathfrak{B}$ can be identified with $(p\otimes \id)(\fA\otimes \K)(p\otimes \id)$.

We first notice that $\cL(T)$, $\cR(T)$ and $\Li(T)$ are separable, thus containing a positive element. Both $\cL(T)$ and $\cR(T)$ are full corners of $\Li(T)$. Let $p$ be the projection such that $p\Li(T)p=\cL(T)$ and $v$ the tripotent from above with $vv^*=p\otimes\id$ and $v^*v=\id$. We can identify $\cL(T)\otimes \K$ with $(p\otimes\id)\Li(T)(p\otimes\id)$. Let $$\phi:\Li(T)\otimes\K\to\Li(T)\otimes\K,\;\;\;x\mapsto vxv^*$$ and notice that $\im(\phi)\sub \cL(T)\otimes\K$ since $(p\otimes\id)(\phi(x))(p\otimes\id)=vv^*vxv^*vv^*=\phi(x)$ for all $x\in\Li(T)\otimes\K$. The linear mapping $\phi$ is bijective with inverse $\phi^{-1}:\Li(T)\otimes\K\to\Li(T)\otimes\K$, $x\mapsto v^*xv$ and $\phi$ is an algebra homomorphism since $\phi(xy)=vxyv^*=vxv^*vyv^*=\phi(x)\phi(y)$ for all $x,y\in\Li(T)\otimes\K$. Since $K_0$ is stable this gives us an isomorphism from $K_0(\cL(T))$ to $K_0(\Li(T))$. The isomorphism from $K_0(\cR(T))$ to $K_0(\Li(T))$ is constructed in a similar way.
 \end{proof}
\end{prop}

We discovered different ways to prove the above lemma. Using the theory of stably isomorphic $C^*$-algebras developed by Brown, as used above, was the proof which needed the least mathematical machinery. The two others involve the use of $KK$-theory, one uses even ideas worked out in the context of $KK$-theory for Banach algebras developed in \cite{Paravicini-MoritaequivalencesandKKtheoryforBanachalgebras}. Both build on the fact that $T$ can be interpreted as a Morita equivalence between $\cL(T)$ and $\cR(T)$, and thus $T$ induces an isomorphism in $KK$-theory.

\begin{cor}\label{kanonische Einbettungen sin K_0-isos}
Let $T$ be a separable TRO. The canonical embeddings $\iota_\cL:\cL(T)\to\Li(T)$ and $\iota_\cR:\cR(T)\to\Li(T)$ induce isomorphisms $$K_0(\iota_\cL):K_0(T)\to K_0(\Li(T))$$ and $$K_0(\iota_\cR):K_0(\cR(T))\to K_0(\Li(T)).$$
\begin{proof}
Choose a system $\{e_{i,j}:i,j\in\N\}$ of matrix-units of $\K$ and consider the commutative diagram \[
\text{\xymatrix{ \cL(T)\ar[d]\ar[rr]^{\iota_\cL}
 &&\Li(T) \ar[d]
\\
\cL(T)\otimes \K\ar[rr]_{\iota_\cL\otimes \id}
 &&\Li(T)\otimes \K}}
\] where the vertical maps are given by $a\mapsto a\otimes e_{1,1}$. On the level of $K_0$, the vertical maps become isomorphisms (this is just the stability of $K_0$) so we only have to prove that $K_0(\iota_\cL\otimes \id)$ is an isomorphism, but this follows from the proof of \ref{moritaisosatz} since $\iota_\cL\otimes\id$ induces the same map on the $K_0$-level as $\phi^{-1}$, by \cite{Brown-StableisomorphismofhereditarysubalgebrasofCalgebras}, Lemma 2.7.
\end{proof}
\end{cor}

\begin{defi}
For a separable TRO $T$ the isomorphism \[\eta_T:=K_0(\iota_{\mathcal{L}})^{-1}\circ K_0(\iota_{\mathcal{R}}):K_0(\mathcal{R}(T))\to K_0(T)\]  is said to be the \textbf{Morita isomorphism}\index{Morita isomorphism} of the left and right $K_0$-groups of $T$.
\end{defi}

The next Lemma shows that the Morita isomorphism respects the group homomorphisms induced by TRO-homomorphisms. This naturality becomes important in the next chapter.

\begin{lemma}\label{natural morita iso}
Let $\vp: T\to U$ be a TRO-homomorphism of separable TROs, then the diagram
\[
\text{\xymatrix{ K_0(\mathcal{R}(T))\ar[d]_{K_0(\mathcal{R}(\vp))}\ar[rr]^{\eta_T}
 &&K_0(T) \ar[d]^{K_0(\vp)}
\\
K_0(\mathcal{R}(U))\ar[rr]_{\eta_U}
 && K_0(U)}}
\]
commutes.
\begin{proof}
This is just part of the commuting diagram \[\text{\xymatrix{ \mathcal{R}(T)\ar[d]_{\mathcal{R}(\vp)}\ar[rr]^{\iota_{\mathcal{R}(T)}}
 &&\Li(T)\ar[d]^{\Li(\vp)}&&\mathcal{L}(T)\ar[ll]_{\iota_{\mathcal{L}(T)}}  \ar[d]^{\mathcal{L}(\vp)}
\\
\mathcal{R}(U)\ar[rr]_{\iota_{\mathcal{R}(U)}}
&&\Li(U)  && \mathcal{L}(U)\ar[ll]^{\iota_{\mathcal{L}(U)}}}}
\] under the functor $K_0$, where the horizontal arrows become isomorphisms due to Corollary \ref{kanonische Einbettungen sin K_0-isos}.
\end{proof}
\end{lemma}

We introduce another notion which is a ternary generalization of a concept which has proven useful in $C^*$-theory.

In $C^*$-algebra theory the scale, defined as the set containing all Murray-von Neumann equivalence classes of projections in the original $C^*$-algebra, is a tool to control the dimension of the $C^*$-algebra (at least in some important cases). We introduce the notion of a double-scale, since a single scale is not sufficient for our intent. A very good example is the TRO $\M_{n,m}$, determined not only by one dimension but by the pair $(n,m)\in\N\times\N$ or equivalently by the dimensions of $\cL(\M_{n,m})=\M_n$ and $\cR(\M_{n,m})=\M_m$. We therefore consider the scales in the left and right $C^*$-algebra simultaneously, transporting the scale of the right $C^*$-algebra with the Morita isomorphism to the $K_0$-group of the TRO.
\begin{defi}
For a TRO  $T$ the \textbf{positive cone}\index{positive cone} of $K_0(T)$ is the set
\[K_0(T)_+:=\left\{[p]\in K_0(T):p\in \bigcup_{n\in\N}M_n(\mathcal{L}(T)), \;p \text{ is a projection}\right\}.\]
\end{defi}

\begin{defi}
Let $T$ be a separable TRO with $K_0$-group $K_0(T)$ and positive cone $K_0(T)_+$.
The \textbf{left scale}\index{left scale} of $K_0(T)$ is the set
\begin{equation*}
\Sigma^\mathcal{L}(T):=\left\{[p]\in K_0(T):p\text{ is a projection in } \mathcal{L}(T)\right\}\sub K_0(T)_+.
\end{equation*}
Let $\eta_T$ be the Morita isomorphism of $T$. We define the \textbf{right scale}\index{right scale} of $K_0(T)$ to be
\begin{equation*}
\Sigma^\mathcal{R}(T):=\eta_T\left(\left\{[p]\in K_0(\mathcal{R}(T)):p\text{ is a projection in } \mathcal{R}(T)\right\}\right)\sub K_0(T)_+.
\end{equation*}
The quadruple \[
\left(K_0(T),K_0(T)_+,\Sigma^\mathcal{L}(T),\Sigma^\mathcal{R}(T)\right)
\] is called \textbf{double-scaled ordered $\mathbf{K_0}$-group}\index{double-scaled ordered $K_0$-group} of $T$.

The homomorphisms of double-scaled ordered $K_0$-groups are those positive group homomorphisms which map the left scale into the left scale and the right scale into the right scale.
\end{defi}

If $\vp:T\to U$ is a TRO-homomorphism $K_0(\vp):K_0(T)\to K_0(U)$ is positive with $K_0(\vp)(\Sigma^\cL(T))\sub \Sigma^\cL(U)$, since $\cL(\vp)$ maps projections to projections. In addition $K_0(\vp)(\Sigma^\cR(T))\sub \Sigma^\cR(U)$ because $K_0(\vp)\circ\eta_T=\eta_U\circ K_0(\cR(\vp))$ holds by Lemma \ref{natural morita iso}.

\begin{exa}\label{endl dim TROs und K_0}
\emph{
Let $U$ be a finite dimensional TRO. By \cite{Smith-Finitedimensionalinjectiveoperatorspaces} we can assume that there exists a $k\in\N$ such that \[
U=\M_{n_1,m_1}\oplus\ldots\oplus\M_{n_k,m_k}.
\] We get with the above and Proposition \ref{K_0 funktor eigenschaften} that
\[K_0(U)=\Z^k\;\text{ and }\;K_0(U)_+=\N^k_0,\]
 \[
\Sigma^\cL(U)=\{(\alpha_1,\ldots,\alpha_k)\in\N^k_0:\alpha_i\leq n_i\text{ for } 1\leq i\leq k \}
\]
and
\[
\Sigma^\cR(U)=\{(\beta_1,\ldots,\beta_k)\in\N^k_0:\beta_i\leq m_i\text{ for } 1\leq i\leq k\}.
\]}
\end{exa}

The next Proposition can proved very similar to $C^*$-theory. We refer the interested reader to \cite{Bohle}.

\begin{prop}\label{K0 TRO homomorphismen}
Suppose $$\varphi: T=\bigoplus_{i=1}^p \mathbb{M}_{n_i,m_i}\to
U=\bigoplus_{j=1}^q \mathbb{M}_{l_j,k_j}$$ is a TRO-homomorphism and let $$\varphi_j:\bigoplus_{i=1}^p \mathbb{M}_{n_i,m_i}\rightarrow
\mathbb{M}_{l_j,k_j}$$  be the restrictions of $\varphi$ for $j=1,\ldots,q$ with $\varphi=\varphi_1+\ldots+\varphi_q$. For every $i\in\{1,\ldots,p\}$ let $\iota_i$ be the embedding of $\mathbb{M}_{n_i,m_i}$ into $\bigoplus_{i=1}^p
\mathbb{M}_{n_i,m_i}$ and
$\varphi_{i,j}:=\varphi_j\circ\iota_i$. The induced group homomorphism $K_0(\varphi):\Z^p\to\Z^q$ is given by the $q\times p$ matrix
\[K_0(\vp)=(\alpha_{i,j})_{i,j},\] where
$\alpha_{i,j}$ is the multiplicity of $\varphi_{i,j}$ for $1\leq i\leq p$, $1\leq j\leq q$.
\end{prop}

\begin{prop}\label{lifting of TRO homs}
Let $T=\bigoplus_{i=1}^p \mathbb{M}_{n_i,m_i}$ and
$U=\bigoplus_{j=1}^q\mathbb{M}_{k_j,l_j}$ be finite dimensional ternary rings of operators and
{\allowdisplaybreaks
$$\left(K_0(T),K_0(T)_+,\Sigma^\mathcal{L}(T),\Sigma^\mathcal{R}(T)\right)=
\left(\Z^p,\N_0^p,\prod_{i=1}^p\{0,\ldots,n_i\},\prod_{i=1}^p\{0,\ldots,m_i\}\right)$$}
 and
$$\left(K_0(U),K_0(U)_+,\Sigma^\mathcal{L}(U),\Sigma^\mathcal{R}(U)\right)
=\left(\Z^q,\N_0^q,\prod_{j=1}^q\{0,\ldots,k_j\},\prod_{j=1}^q\{0,\ldots,l_j\}\right)$$ their double-scaled ordered $K_0$-groups.
Let $\alpha:K_0(T)\rightarrow K_0(U)$ be a homomorphism of double-scaled groups.
\begin{description}
  \item[(a)] The homomorphism $\alpha$ can be represented as a $q\times p$-matrix $(a_{i,j})_{i,j}$ with entries $a_{i,j}\in\N_0$.
  \item[(b)] For all $(z_1,\ldots z_p)\in \Sigma^\mathcal{L}(T)$ we have $\sum_{j=1}^pa_{i,j}z_j\leq
  k_i$ for all $i=1,\ldots,q$.
  \item[(c)] For all $(z_1,\ldots z_p)\in \Sigma^\mathcal{R}(T)$ we have $\sum_{j=1}^pa_{i,j}z_j\leq
  l_i$ for all $i=1,\ldots,q$.
  \item[(d)] There exists a TRO-homomorphism $\varphi:T\to U$ with
  $K_0(\varphi)=(a_{i,j})_{i,j}$.
\end{description}
\begin{proof}
All group homomorphisms from $\mathbb{Z}^p$ to $\mathbb{Z}^q$ can be viewed as
$q\times p$-matrices with entries in $\mathbb{Z}$. The homomorphism $\alpha$ is a homomorphism of ordered groups and therefore maps $K_0(T)_+=\N_0^
p$ to $K_0(U)_+=\N_0^
q$, thus all entries in the matrix have to be positive or $0$.

To prove (b) let $x=(z_1,\ldots z_p)\in \Sigma^\mathcal{L}(T)$. Since $\alpha
(\Sigma^\mathcal{L}(T))\subseteq \Sigma^\mathcal{L}(U)$ holds, we see that
$$\alpha(x)=\left(\sum a_{1,i}z_i,\ldots,\sum
a_{q,i}z_i\right)\leq(k_1,\ldots,k_q).$$

An analogous argument shows (c).

For the proof of (d) let $\varphi$ be the direct sum $\varphi:=\varphi_1\oplus\ldots\oplus\varphi_q$, where
$\varphi_j:T\rightarrow \mathbb{M}_{k_j,l_j}$ is defined via
\[\varphi_j(x_1\oplus\ldots \oplus x_p):=\diago(\underbrace{x_1,\ldots,x_1}_{a_{1,j} times},\ldots,\underbrace{x_p,\ldots x_p}_{a_{p,j} times},0,\ldots,0),\]
for $j=1,\ldots,q$.
These TRO-homomorphisms are well-defined by (b) and (c) and using Proposition \ref{K0 TRO homomorphismen}
 we get $K_0(\varphi)=\alpha$.
\end{proof}
\end{prop}

\begin{prop}\label{isomorphe endl dim TROs}
Two finite dimensional TROs are isomorphic if and only if their double-scaled ordered groups are isomorphic.
\begin{proof}
Let $T$ and $U$ be finite dimensional TROs. If $\varphi:T\to U$ is a TRO-isomorphism, then $K_0(\vp)$ becomes an isomorphism of double-scaled ordered groups.

 If on the other hand $\left(K_0(T),K_0(T)_+,\Sigma^\cL(T),\Sigma^\cR(T)\right)$ is the scaled ordered group of $T$, then we know from Example \ref{endl dim TROs und K_0} that there exist natural numbers $k,n_1,\ldots,n_k,m_1,\ldots m_k$ such that $K_0(T)\simeq\Z^k$, $K_0(T)_+\simeq\N_0^k$, $\Sigma^\cL(T)\simeq\prod_{i=1}^k\{0,\ldots,n_i\}$ and $\Sigma^\cR(T)\simeq\prod_{i=1}^k\{0,\ldots,m_i\}$. Since every finite dimensional TRO is isomorphic to the direct sum of rectangular matrix algebras, the double-scaled ordered $K_0$-group of $T$ carries all the necessary data to recover $T$ up to isomorphism. If $T\simeq\bigoplus_{i=1}^l\M_{r_i,s_i}$, then $l=k$ and (after maybe changing the summation order) $(r_i,s_i)=(m_i,n_i)$ for $i=1,\ldots,k$.
\end{proof}
\end{prop}

In \cite{Bohle} we extended the previous result to inductive limits of finite dimensional TROs.

\begin{rem}
\emph{One can easily get the impression that ternary $K$-theory is simply the doubled $K$-theory of two Morita equivalent $C^*$-algebras and not so much an invariant of the TRO itself. But this is not the case as the following simple but illuminating twofold example shows.
The TROs
\begin{align*}
  T=\M_{1,2}\oplus\M_{2,1}\;\;\;\;\;\;\text{ and }\;\;\;\;\;\; U=\M_{1,1}\oplus\M_{2,2}
\end{align*}
are non-isomorphic (not even linear isomorphic) with $$\cL(T)=\M_{1}\oplus\M_2=\cL(U)$$
and $$\cR(T)=\M_2\oplus\M_1\simeq\M_1\oplus\M_2=\cR(U).$$
This yields the two non-isomorphic double-scaled ordered groups
\begin{align*}K_0(T)=(\Z^2,\N_0^2,&\{(0,0),(0,1),(0,2),(1,1),(1,2)\},\\
 &\{(0,0),(1,0),(1,1),(2,0),(2,1)\})
\end{align*}
and
\begin{align*}K_0(U)=(\Z^2,\N_0^2,&\{(0,0),(0,1),(0,2),(1,1),(1,2)\},\\
 &\{(0,0),(0,1),(0,2),(1,1),(1,2)\}).\end{align*}}
\emph{
The example shows that the scaled $K$-theory of TROs is not a fusion of the scaled $K$-theories of two Morita equivalent $C^*$-algebras, but can distinguish between different TROs, even if they have isomorphic left and right $C^*$-algebras.}
\end{rem}

\section{K-theory for JB*-triple systems}

In \cite{BoWe1} we associated to every $JB^*$-triple system $Z$ an, up to TRO-isomorphism, unique pair $(T^*(Z),\rho_Z)$, where $T^*(Z)$ is a TRO and $\rho_Z:Z\to T^*(Z)$ is a $JB^*$-triple homomorphism with the following two universal properties: (i) for every $JB^*$-triple homomorphism $\vp:Z\to T$ to a TRO $T$ there exists a TRO-homomorphism $T^*(\vp):T^*(Z)\to T$ with $T^*(\vp)\circ\rho_Z=\vp$; (ii) $\rho_Z(Z)$ generates $T^*(Z)$ as a TRO.

The universal properties of $T^*$ yield a functor $\tau$ from the category of $JB^*$-triple homomorphisms to the category of TROs, if we map a $JB^*$-triple system $Z$ to its
universal enveloping TRO $\tau(Z):=T^*(Z)$ and $JB^*$-triple
homomorphisms $\varphi:Z_1\to Z_2$ to TRO-homomorphisms
$\tau(\varphi):=T^*(\rho_{Z_2}\circ\varphi):\tau(Z_1)\to \tau(Z_2)$.
 We examine the functorial properties of the mapping $\tau$, or more exactly its restriction to the category of $JC^*$-triple systems. As it turns out $\tau$ is homotopy invariant, continuous, additive and exact.

We use the results of the previous two chapters to define the $K$-groups of $JB^*$-triple systems. For a given $JB^*$-triple system $Z$ we define the $i$th $K$-group of $Z$ to be the $i$th (ternary) $K$-group of its universal enveloping TRO $T^*(Z)$. By the results of the first subsection we obtain covariant, continuous, half-exact, split-exact and homotopy invariant functors on the subcategory of $JC^*$-triple systems.

Using the theory of grids combined with our ordered $K$-theory for TROs we are able to define an invariant for atomic $JBW^*$-triple systems given by a tuple $$
\left(
K_0^{\text{\tiny JB*}}(Z),K_0^{\text{\tiny JB*}}(Z)_+,\Sigma^{\text{\tiny JB*}}
_\mathcal{L}(Z),\Sigma^{\text{\tiny JB*}}_\mathcal{R}(Z),\Gamma(Z)
\right),
$$
where $\left(K_0^{\text{\tiny JB*}}(Z),K_0^{\text{\tiny JB*}}(Z)_+,\Sigma^{\text{\tiny JB*}}
_\mathcal{L}(Z),\Sigma^{\text{\tiny JB*}}_\mathcal{R}(Z)\right)$ is the double-scaled ordered $K_0$-group of $T^*(Z)$ and $\Gamma(Z)$ is the set of equivalence classes in $K_0^{\text{\tiny JB*}}(Z)_+$ which stem from a grid spanning $Z$. We show that this invariant is a well defined, complete isomorphism invariant for finite dimensional $JC^*$-triple systems. We prove this by computing the invariants of all finite dimensional $JC^*$-triple systems.

\subsection{Functorial properties}\label{section functorial properties}

The next Lemma is needed to prove exactness of $\tau$.

\begin{lemma}\label{TROideallemma}
Let $Z$ be a $JC^*$-triple system and $T$ a TRO such that $Z$
generates $T$ as a TRO. If $I\subseteq Z$ is a $JB^*$-triple ideal,
then the TRO $[I]$ generated by $I$ in $T$ is a TRO-ideal in $T$.

\begin{proof}
We have to show that $[I]T^*T+T[I]^*T+TT^*[I]\subseteq [I]$. Since
$I$ generates $[I]$ and $Z$ generates $T$ it suffices to show that
for all $i\in I$ and for all $x,y\in Z$: $ix^*y$, $xi^*y,
xy^*i\in[I].$ To prove this we first show that for all $i_1,i_2\in
I$, $z\in Z$: $i_1 i_2^*z,$ $zi_2^*i_1\in[I].$ There exists  $j\in I$ with
$\{j,j,j\}=jj^*j=i_2$. Thus
\begin{align*}
i_1i_2^*z=i_1j^*jj^*z &=  2i_1j^*\{j,j,z\}-i_1j^*zj^*j\displaybreak[0]\\
&=2i_1j^*\{j,j,z\}-i_1\{j,z,j\}^*j\in[I]
\end{align*}
and similarly $zi_2^*i_1\in [I]$.
Now let $x,y\in Z$, $i\in I$ and, as above, $j\in I$ with $jj^*j=i$.
We have
\begin{align*}
xi^*y=xj^*jj^*y&=2\{x,j,j\}j^*y-jj^*xj^*y\displaybreak[0]\\
&=2\{x,j,j\}j^*y - j\{j,x,j\}^*y\in[I],
\end{align*}
and similarly we get $ix^*y\in [I]$ and $xy^*i\in [I]$.
\end{proof}
\end{lemma}

Likewise to the case of $JC$-algebras (cf.\
\cite{HancheOlsen-OnthestructureandtensorproductsofJCalgebras}) the
functor $\tau$ is exact.

\begin{thm}\label{exactheit von T*}
Every exact sequence of $JC^*$-triple systems\[ 0\to I\to Z\to
Z/I\to 0,\] where I is a $JB^*$-triple ideal of $Z$, induces an
exact sequence of the corresponding universal enveloping TROs:\[0\to
\tau(I)\to \tau(Z)\to \tau(Z/I)\to 0.\]

\begin{proof}
Let $\iota: I\to Z$ and $\pi: Z\to Z/I$ be the canonical injection
and quotient homomorphism, respectively.

We first show exactness at $\tau(Z/I)$: The TRO $\tau(Z)$ is
generated by $\rho_Z(Z)$ and the TRO $\tau(Z/I)$ by
\begin{align*}
\rho_{Z/I}(Z/I)&=\rho_{Z/I}(\pi(Z))\displaybreak[0]\\
&=T^*(\rho_{Z/I}\circ\pi)(\rho_Z(Z))\displaybreak[0]\\
&=\tau(\pi)(\rho_Z(Z)).
\end{align*}
Next we show exactness at $\tau(Z)$: We have $\tau(\pi)\circ
\tau(\iota)=\tau(\pi\circ \iota)=0$ by functoriality. Let
$\widetilde{I}:=\tau(\iota)(\tau(I))$, then $\widetilde{I}$ is a TRO and
the $JB^*$-triple ideal $\rho_Z(\iota(I))\subseteq\rho_Z(Z)$
generates $\widetilde{I}$ as a TRO. By Lemma \ref{TROideallemma} the subTRO
$\widetilde{I}$ is a TRO-ideal of $\tau(Z)$.

Let $\widetilde{\pi}:\tau(Z)\to \tau(Z)/\widetilde{I}$ be the quotient
homomorphism onto the TRO
$\tau(Z)/\widetilde{I}$, then $\widetilde{\pi}\circ \rho_Z\circ
\iota=0$, since $\widetilde{I}$ is generated by $\rho_Z(\iota (I))$.
Therefore the $JB^*$-triple homomorphism
$\widetilde{\pi}\circ\rho_Z$ induces a $JB^*$-triple homomorphism
$\varphi:Z/I\to \tau(Z)/\widetilde{I}$, which induces the
TRO-homomorphism $T^*(\varphi):\tau(Z/I)\to \tau(Z)/\widetilde{I}$.
In the diagram
\[\xymatrix{
Z \ar[rr]^-{\rho_Z}\ar[rrdd]_{\rho_{Z/I}\circ\pi}&&
\tau(Z)\ar[rr]^-{\widetilde{\pi}}\ar[dd]^{\tau(\pi)} &&
\tau(Z)/\widetilde{I}
\\\\
& &\tau(Z/I)\ar[rruu]_-{T^*(\varphi)}}\] the left triangle trivially commutes and the outer triangle commutes by the definition of $\vp$. If $z\in\rho_Z(Z)\sub\tau(Z)$, we find an element $x\in Z$ with $\rho_Z(x)=z$. Since the left and the outer triangle commute we have \begin{align*}
  \widetilde{\pi}(z)&=\widetilde{\pi}(\rho_Z(x))\\
  &=T^*(\vp)\circ\rho_{Z/I}\circ\pi(x)\\
  &=T^*(\vp)\circ\tau(\pi)(\rho_Z(x))\\
  &=T^*(\vp)\circ\tau(\pi)(z).
\end{align*}
Now $\rho_Z(Z)$ generates $\tau(Z)$ and thus the right triangle commutes. We obtain the desired inclusion
$\ker(\tau(\pi))\subseteq\widetilde{I}=\tau(\iota)(\tau(I))$.

Finally we show exactness at $\tau(I)$: Let $H$ be a Hilbert space
and $\alpha:\tau(I)\to B(H)$ an injective TRO-homomorphism. Then
$\alpha\circ\rho_I:I\to B(H)$ is an injective $JB^*$-homomorphism
and we get by \cite{BoWe1} Lemma 2.5 an extension
$\overline{\alpha}:Z\to B(H)$ of $\alpha\circ\rho_I$. This
$JB^*$-homomorphism lifts to a TRO-homomorphism
$T^*(\overline{\alpha}):\tau(Z)\to B(H)$. With the universal
property we get $T^*(\overline{\alpha})\circ \tau(\iota)=\alpha$.
Since $\alpha$ is injective, so is $\tau(\iota):\tau(I)\to \tau(Z)$.
\end{proof}
\end{thm}

\begin{defi}
Let $Z_1$ and $Z_2$ be $JB^*$-triple systems and
$\alpha,\beta:Z_1\to Z_2$ $JB^*$-triple homomorphisms. The mappings
$\alpha$ and $\beta$ are called \textbf{JB*-homotopic}\index{homotopic $JB^*$-triple homomorphisms}, denoted
$\alpha\sim_h^{JB^*}\beta$, when there is a path
$(\gamma_t)_{t\in[0,1]}$ of $JB^*$-triple homomorphisms
$\gamma_t:Z_1\to Z_2$ such that $t\mapsto \gamma_t(z)$ is a norm
continuous path in $Z_2$ for every $z\in Z_1$ with
$\gamma_0=\alpha$, $\gamma_1=\beta$.

A $JB^*$-triple homomorphism $\alpha:Z_1\to Z_2$ is called a
\textbf{JB*-homotopy equivalence} when there is a $JB^*$-triple homomorphism
$\beta:Z_2\to Z_1$ such that $\alpha\circ \beta$ and $\beta\circ
\alpha$ both are homotopic to the identity.
\end{defi}

\begin{prop}
Let $\alpha,\beta:Z_1\to Z_2$ be homotopic $JB^*$-triple
homomorphisms. If $(T^*(Z_1),\rho_{Z_1})$ and
$(T^*(Z_2),\rho_{Z_2})$ are the corresponding universal enveloping
TROs, then the functor $\tau$ induces a TRO homotopy between
$\tau(\alpha)$ and $\tau(\beta)$.
\begin{proof}
Let $(\gamma_t)_{t\in[0,1]}$ be a pointwise continuous path in $Z_2$ which
connects $\alpha$ and $\beta$ and
$(\widetilde{\gamma}_t)_{t\in[0,1]}$ be the path defined by
$\widetilde{\gamma}_t:=\tau(\gamma_t):\tau(Z_1)\to \tau(Z_2)$.
Obviously $(\widetilde{\gamma}_t)$ connects $\tau(\alpha)$ and
$\tau(\beta)$, so the only thing to show is that
$t\mapsto\widetilde{\gamma}_t(z)$ defines
a norm continuous path in $\tau(Z_2)$ for every $z\in \tau(Z_1)$. Since $\tau(Z_1)$ is
generated by $\rho_{Z_1}(Z_1)$ we can assume w.l.o.g.\ that
$$z=\rho_{Z_1}(z_1)\rho_{Z_1}(z_2)^*\rho_{Z_1}(z_3)\ldots\rho_{Z_1}(z_{2n})^*\rho_{Z_1}(z_{2n+1})$$
with $z_i\in Z_1$. Then
\begin{align*}
\widetilde{\gamma}_t(z) &=
\widetilde{\gamma}_t(\rho_{Z_1}(z_1))\widetilde{\gamma}_t(\rho_{Z_1}(z_2))^*\widetilde{\gamma}_t(\rho_{Z_1}(z_3))\ldots\widetilde{\gamma}_t(\rho_{Z_1}(z_{2n}))^*\widetilde{\gamma}_t(\rho_{Z_1}(z_{2n+1}))\\
&=
\rho_{Z_2}(\gamma_t(z_1))\rho_{Z_2}(\gamma_t(z_2))^*\rho_{Z_2}(\gamma_t(z_3))\ldots\rho_{Z_2}(\gamma_t(z_{2n}))^*\rho_{Z_2}(\gamma_t(z_{2n+1})),
\end{align*}
which is norm continuous in t.
\end{proof}
\end{prop}

With the help of our functor $\tau$ we are able to show that inductive limits exist in the category of $JC^*$-triple systems.

\begin{lemma}
Let $((Z_n),(\vp_n))$ be an inductive sequence in the category of $JC^*$-triple systems, then $((\tau(Z_n)),(\tau(\vp_n)))$ is an inductive sequence of TROs. Moreover, if $(T_\infty,(\mu_n))$ is the inductive limit of $((\tau(Z_n)),(\tau(\vp_n)))$ in the category of TROs, then $(Z_\infty,(\nu_n))$ is the inductive limit in the category of $JC^*$-triple systems, where $$Z_\infty:=\overline{\bigcup_{n=1}^{\infty}\mu_n(\rho_{Z_n}(Z_n))}$$  with homomorphisms $\nu_n:=\mu_n\circ\rho_{Z_n}:Z_n\to Z_\infty$ for all $n\in\N$.
\begin{proof}
It is straightforward to check that $((\tau(Z_n)),(\tau(\vp_n)))$ is an inductive sequence of TROs. This sequence of TROs has an inductive limit by Proposition \ref{inductive lim fuer TROs} and an argument similar to its proof shows the rest of the lemma.
\end{proof}
\end{lemma}

Now that we have proved the existence of inductive limits in the category of $JC^*$-triples we can prove the continuity of the functor $\tau$.

\begin{prop}
Let $((Z_n),(\vp_n))$ be an inductive sequence in the category of $JC^*$-triple systems. If the pair $(Z_\infty,(\mu_n))$ is the inductive limit of $((Z_n),(\vp_n))$, then $(\tau(Z_\infty),(\tau(\mu_n)))$ is the inductive limit of the induced sequence of TROs.
\begin{proof}
We know by functoriality of $\tau$ that $$\tau(\mu_{n+1})\circ\tau(\vp_n)=\tau(\mu_{n+1}\circ\vp_n)=\tau(\mu_n)$$
for all $n\in\N$, so \begin{equation*}
\text{\xymatrix{ \tau(Z_n)\ar[dr]_{\tau(\mu_n)}
\ar[rr]^-{\tau(\vp_n)} && \tau(Z_{n+1})\ar[dl]^{\tau(\mu_{n+1})}
\\
 &\tau(Z_{\infty}) }}
 \end{equation*}
 commutes for all $n\in\N$.

  Let $(T_\infty,(\lambda_n))$ be another TRO system with $\lambda_n=\lambda_{n+1}\circ\tau(\vp_n)$ for all $n\in\N$. We
 have to show that there exists a unique TRO-homomorphism $\lambda$ making
\begin{equation}\label{induktiver limes b}
\text{\xymatrix{
 &\tau(Z_n)\ar[dr]^{\lambda_n}\ar[dl]_{\tau(\mu_{n})}
\\
\tau(Z_{\infty})
\ar[rr]^-{\lambda}
&&
T_\infty
 }}\end{equation} commute for all $n\in\N$. We first notice that the commutative diagram

\begin{equation*}
\text{\xymatrix{
 Z_n\ar[rr]^{\mu_n}\ar[d]_{\rho_{Z_{n}}}&& Z_{n+1}\ar[d]^{\rho_{Z_{n+1}}}\\
 T^*(Z_n)\ar[rr]^{\tau(\mu_n)}\ar[dr]_{\lambda_n}&& T^*(Z_{n+1})\ar[dl]^{\lambda_{n+1}}\\
 &T_\infty
 }}\end{equation*}
induces by the universal property of the inductive limit a unique $JB^*$-triple homomorphism $\lambda:Z_\infty\to T_\infty$ such that

\begin{equation*}
\text{\xymatrix{
 Z_n\ar[rr]^{\rho_{Z_n}}\ar[d]_{\mu_n}&& T^*(Z_n)\ar[d]^{\lambda_n}\\
 Z_\infty\ar[rr]^{\lambda}&& T_\infty
 }}\end{equation*}
commutes for all $n\in\N$. This induces a unique TRO-homomorphism $T^*(\lambda):T^*(Z_\infty)\to T_\infty$ such that (\ref{induktiver limes b})
 commutes for all $n\in\N$, which shows that $(\tau(Z_\infty),(\tau(\mu_n)))$ is the inductive limit of $(\tau(Z_n),(\tau(\vp_n)))$.
\end{proof}
\end{prop}

\subsection{The K-groups of a JB*-triple system}

\begin{defi}
Let $Z$ be a $JB^*$-triple system and $T^*(Z)$ its universal enveloping TRO. We define the $\mathbf{i}$\textbf{th K-group}\index{$K_0$-group of a $JB^*$-triple system}, $i\in \N_0$, of $Z$ by \[K_i^{\text{\tiny JB*}}(Z):=K_i(\tau(Z)).\]
If $\vp:Z\to W$ is a $JB^*$-triple homomorphism, then $\vp$ induces for all $i\in\N_0$ group homomorphisms $$
K_i^{\text{\tiny JB*}}(\vp):K_i^{\text{\tiny JB*}}(T)\to K_i^{\text{\tiny JB*}}(W),
$$ defined by $K_i^{\text{\tiny JB*}}(\vp)=K_i(\tau(\vp))$.
\end{defi}

The next proposition follows immediately from Section \ref{section functorial properties} and Proposition \ref{K_0 funktor eigenschaften}.

\begin{prop}
Let $i\in\N_0$, then
$K_i^{\text{\tiny JB*}}$ is a covariant functor from the category of $JB^*$-triple systems to the category of Abelian groups.
 Restricted to the category of $JC^*$-triple systems the functors $K_i^{\text{\tiny JB*}}$ are half-exact, split-exact, homotopy invariant, continuous and Bott periodic.
 \end{prop}

\section{A complete isomorphism invariant}

A $JB^*$-triple $Z$ system which is a dual Banach space is called a
$\mathbf{JBW^*}$-\textbf{triple system}\index{$JBW^*$-triple
system}.

Let $Z$ be a $JB^*$-triple system. A tripotent $e\in Z$ is called \textbf{minimal}\index{minimal tripotent}
if \[\{e,Z,e\}=\mathbb{C}e.\]
A $JBW^*$-triple system is called \textbf{atomic}\index{atomic $JBW^*$-triple system} if it is the $w^*$-closed linear span of its minimal tripotents. The importance of atomic $JBW^*$-triples comes from the fact, that they have a decomposition $Z=\bigoplus_i Z_i$ into a $l^\infty$-direct sum of (possibly infinite dimensional) Cartan factors and each of those Cartan factors is spanned by a standard grid (cf. \cite{Neher-Jordantriplesystemsbythegridapproach}). We will give examples of the standard grids of the finite dimensional Cartan factors in Chapter 5.

\begin{defi}
Let $Z$ be an atomic $JBW^*$-triple system spanned by a grid $\cG$. The \textbf{K-grid invariant}\index{K-grid invariant of a $JB^*$-triple system} of $Z$ is the tuple \[\mathcal{KG}(Z):=\left(K_0^{\text{\tiny JB*}}(Z),K_0^{\text{\tiny JB*}}(Z)_+,\Sigma^{\text{\tiny JB*}}
_\mathcal{L}(Z),\Sigma^{\text{\tiny JB*}}_\mathcal{R}(Z),\Gamma(Z)\right),\] where $K_0^{\text{\tiny JB*}}(Z)_+:=K_0(T^*(Z))_+$, $\Sigma_\cL^{\text{\tiny JB*}}(Z)$ and $\Sigma^{\text{\tiny JB*}}_\cR(Z)$ are the left and right scale of the TRO $T^*(Z)$ and $\Gamma(Z)$ is the set of equivalence classes \[
\Gamma(Z):=\left\{[\rho_Z(g)\rho_Z(g)^*]\in K_0^{\text{\tiny JB*}}(Z): g\in \cG\right\}\sub \Sigma_\mathcal{L}^{\text{\tiny JB*}}(Z).
\]

Let $\varphi:K_0^{\text{\tiny JB*}}(Z_1)\to K_0^{\text{\tiny JB*}}(Z_2)$ be a group homomorphism. We say that $\varphi$ is a \textbf{K-grid} \textbf{isomorphism of} $\mathbf{K_0}$ \textbf{-groups} \index{K-grid isomorphism of $K_0$-groups} if $\varphi$ is a group isomorphism with $\vp(K_0^{\text{\tiny JB*}}(Z)_+)=K_0^{\text{\tiny JB*}}(W)_+$, $\varphi(\Sigma^{\text{\tiny JB*}}_\mathcal{L}(Z_1))=\Sigma^{\text{\tiny JB*}}_\mathcal{L}(Z_2)$, $\varphi(\Sigma^{\text{\tiny JB*}}_\mathcal{R}(Z_1))=\Sigma^{\text{\tiny JB*}}_\mathcal{R}(Z_2)$ and $\varphi(\Gamma(Z_1))=\Gamma(Z_2)$.
\end{defi}

Our choice of the equivalence classes of the grid elements as additional classifying data is of course not by chance. The grids, as shown by Neher (cf.\ for example \cite{Neher-Systemesderacines3gradues}, \cite{Neher-3gradedrootsystemsandgridsinJordantriplesystems} and \cite{NeherErhard-Liealgebrasgradedby3gradedrootsystemsandJordanpairscoveredbygrids}), are the Jordan analogue of the root systems which were used by \'{E}. Cartan to classify the bounded symmetric spaces in finite dimensions.

The notion of the $K$-grid invariant can be extended to general $JB^*$-triple systems: We first recall from \cite{FriedmanRusso-TheGelfandNaimarktheoremforJBtriples}, Theorem D that every $JBW^*$-triple system $Z$ with predual $Z_*$ decomposes into an orthogonal direct sum $Z=A\oplus N$ of $w^*$-closed ideals $A$ and $N$, where $A$ is the $w^*$-closure of the linear span of its minimal tripotents and $N$ does not contain minimal tripotents. The ideal $A$ is called the atomic part of $Z$. Moreover, they showed in the same article that if $\widetilde{A}$ is the atomic part of $Z''$ and if one composes the canonical embedding $\iota:Z\to Z''$ with the canonical projection $\pi:Z''\to \widetilde{A}$, then $\pi\circ\iota$ is a $JB^*$-triple embedding. Therefore we can define the K-grid invariant of a $JB^*$-triple system $Z$.

Let $Z$ be a $JB^*$-triple system. The $\mathbf{JBW^*}$-\textbf{K-grid invariant} of $Z$ is the K-grid invariant of the atomic part of $Z''$.

The definition of the K-grid invariant of a $JBW^*$-triple system and the $JB^*$-K-grid invariant of a $JB^*$-triple system coincide in the case that the triple system is reflexive as a Banach space, especially in finite dimensions, where our main interest lies. It is known (cf.\ \cite{ChuIochum-ComplementationofJordantriplesinvonNeumannalgebras}, Theorem 6) that a $JB^*$-triple system is reflexive if and only if it does not contain a copy of the function space $c_0$.

However we still have to show that the K-grid invariant is well-defined.

\begin{lemma}\label{gridiso}
Let $Z$ be an atomic $JBW^*$-triple system and $\cG_1$ and $\cG_2$ two grids spanning $Z$, then there exists a $JB^*$-triple automorphism of $Z$ mapping $\cG_1$ onto $\cG_2$.
\begin{proof}This follows from the Isomorphism Theorem 3.18 in \cite{Neher-Jordantriplesystemsbythegridapproach}.
\end{proof}
\end{lemma}

\begin{lemma}\label{isomorphiinvarianz fuer ext K_0}
Let $Z_1$ and $Z_2$ be atomic $JBW^*$-triple systems and $\varphi: Z_1\to Z_2$ be a $JB^*$-triple isomorphism. Let $\cG_i$ be a grid spanning $Z_i$, $i=1,2$. Then there exists a $JB^*$-triple isomorphism $\psi:Z_1\to Z_2$ such that the induced map $K_0^{\text{\tiny JB*}}(\psi)$ is a K-grid isomorphism of $K_0$-groups.
\begin{proof}
 We assume w.l.o.g.\ that $Z_1$ and $Z_2$ are factors and therefore $\vp(\cG_1)$ and $\vp(\cG_1)$ are grids of the same type.
 Let $\psi':Z_2\to Z_2$ be the $JB^*$-triple isomorphism from Lemma \ref{gridiso} which maps $\varphi(\cG_1)$ onto $\cG_2$, then $\psi:=\psi'\circ\varphi$ is a $JB^*$-triple isomorphism with $\psi(\cG_1)=\cG_2$. Therefore $K_0^{\text{\tiny JB*}}(\psi):K_0^{\text{\tiny JB*}}(Z_1)\to K_0^{\text{\tiny JB*}}(Z_2)$ is an isomorphism of double-scaled ordered groups which maps $\Gamma(Z_1)$ onto $\Gamma(Z_2)$.
\end{proof}
\end{lemma}

\begin{prop}\label{gammaZistinvariante}
Let $Z_1$ and $Z_2$ be finite dimensional $JC^*$-triple systems. If $\vp:Z_1\to Z_2$ is a $JB^*$-triple isomorphism then $$K_0(\vp)\left(\Gamma(Z_1)\right)=\Gamma(Z_2).$$
\begin{proof}
We can assume that $Z_1$ and $Z_2$ are simple and spanned by grids $\cG_1\sub Z_1$ and $\cG_2\sub Z_2$, which are of the same type. We consider the images of $\cG_1$ and $\cG_2$ in $T^*(Z_2)$, say $\cG_1':=\rho_{Z_2}(\vp(\cG_1))$ and $\cG_2':=\rho_{Z_2}(\cG_2)$, then there exists by Lemma \ref{isomorphiinvarianz fuer ext K_0} a $JB^*$-triple automorphism $\psi$ mapping $\cG_1'$ to $\cG_2'$. By the universal property of the universal enveloping TRO, we know that $\psi$ has to be the restriction of a TRO-automorphism $\tau(\psi)$ of $T^*(Z_2)$. It is known that every TRO-automorphism in finite dimensions is inner, thus there exist unitaries $U$ and $K$ such that $\tau(\psi)(z)=UzK$ for all $z\in T^*(Z_2)$. Especially we have $$
U\cG_1'K=\cG_2'.
$$
\end{proof}
\end{prop}

Thus any isomorphism of finite dimensional $JC^*$-triple systems yields a K-grid isomorphism of their $K_0$-groups, independent of the choice of the grid.

By the direct sum of two K-grid invariants we mean the direct sum of all the components.

\begin{prop}\label{kgrid inv additiv}
Let $Z_1$ and $Z_2$ be atomic $JBW^*$-triple systems, then there exists a K-grid isomorphism of $K_0$-groups
$$\mathcal{KG}(Z_1\oplus Z_2)\simeq \mathcal{KG}(Z_1)\oplus \mathcal{KG}(Z_2).$$

\begin{proof}
We already know that the functor $K_0$ from the category of TROs to the category of Abelian groups is additive (cf.\ Proposition \ref{K_0 funktor eigenschaften}). This also holds for the positive cone and the scales. The functor $\tau$ defined in section \ref{section functorial properties} is additive. If $\cG$ is a grid spanning $Z_1\oplus Z_2$ and $p_i$ is the projection onto $Z_i$, then $p_i\cG$ is a grid which generates $Z_i$ for $i=1,2$. Therefore $\Gamma(Z_1\oplus Z_2)=\Gamma(Z_1)\oplus\Gamma(Z_2)$.
\end{proof}
\end{prop}

\section{Classification}\label{classification of triples}

We determine the K-grid invariants of all finite dimensional $JC^*$-triple systems. We do this by making a case by case study of the K-grid invariants of the finite dimensional Cartan factors of type I--IV. The universal enveloping TROs of the finite dimensional Cartan factors were computed in \cite{BoWe1}.

\subsection{Rectangular factors}
Recall that a finite dimensional rectangular Cartan factor is a $JC^*$-triple system which is isometric to $$C^1_{n,m}=\M_{n,m}$$
for $n,m\in\N$. The standard example of a rectangular grid spanning $\M_{n,m}$ is $$\cG=\{E_{i,j}:1\leq i\leq n,1\leq j\leq m\}.$$

Let $Z$ be $JC^*$-triple system which is isomorphic to the finite dimensional Cartan factor $C^1_{n,m}$. We have to distinguish between the case when $Z$ is a rank $1$ $JB^*$-triple system and the case $1<n,m<\infty$.

\begin{prop} If $Z$ is the finite dimensional type I Cartan factor $Z=C^1_{n,m}$, with $n,m\geq 2$, then $\mathcal{KG}(C^1_{n,m})$ is given by
\[\left(\Z^2,\N_0^2,\{1,\ldots,n\}\times\{1,\ldots,m\},\{1,\ldots,m\}\times\{1,\ldots,n\},\{(1,1)\}\right)
\]
\begin{proof}
We know by \cite{BoWe1} that the universal enveloping TRO of $Z$ is $$T^*(Z)=\mathbb{M}_{n,m}\oplus \mathbb{M}_{m,n}.$$ If we identify $Z$ with the diagonal $$\{(A,A^t):A\in\M_{n,m} \}\sub\mathbb{M}_{n,m}\oplus\M_{m,n},$$ then $Z$ is spanned by the rectangular grid $$\cG=\{(E_{i,j},E_{j,i}):1\leq i\leq n, 1\leq j\leq m\},$$
thus $\Gamma(Z)$ collapses to $$\Gamma(Z)=\{(1,1)\}.$$
The equality
\begin{align*}
\left(K_0^{\text{\tiny JB*}}(Z),K_0^{\text{\tiny JB*}}(Z)_+,\Sigma^{\text{\tiny JB*}}_\cL(Z),\Sigma^{\text{\tiny JB*}}_\cR(Z)\right)=(\Z^2,\N_0^2,&\{1,\ldots,n\}\times\{1,\ldots,m\},\\&\{1,\ldots,m\}\times\{1,\ldots,n\})
\end{align*}
 follows from the K-theory for ternary rings of operators.
\end{proof}
\end{prop}

Recall from \cite{BoWe1} that, if $Z$ is a finite dimensional type I Cartan factor of rank $1$, then its universal enveloping TRO is given by \[T^*(Z)=\bigoplus_{k=1}^n\M_{p_k,q_k},\]
where $p_k=\begin{pmatrix}n\\
k\end{pmatrix}$ and $q_k=\begin{pmatrix}n\\ k-1\end{pmatrix}$ for
$k=1,\ldots,n$. The image of $Z$ under the injection into $T^*(Z)$ is located inside the direct sum of the spaces $H^k_n$,  $$\rho_Z(Z)\sub\bigoplus_{k=1}^n H^k_n\sub\bigoplus_{k=1}^n\M_{p_k,q_k}$$ (see below for details on the spaces $H^k_n$).

\begin{prop}
If $Z$ is isometric to a finite dimensional Hilbert space, i.e.\ $Z=C^1_{1,n}$, $n\in \N$, then $\mathcal{KG}(C^1_{1,n})$ is given by
\begin{align*}
K_0^{\text{\tiny JB*}}(C^1_{1,n})\;\;&=\Z^n,\\
K_0^{\text{\tiny JB*}}(C^1_{1,n})_+&=\N_0^n,\\
\Sigma_\mathcal{L}^{\text{\tiny JB*}}(C^1_{1,n})\;\;&=\prod_{k=1}^n\left\{1,\ldots,\begin{pmatrix}n\\
k\end{pmatrix}\right\},\\
\Sigma_\mathcal{R}^{\text{\tiny JB*}}(C^1_{1,n})\;\;&=\prod_{k=1}^n\left\{1,\ldots,\begin{pmatrix}n\\ k-1\end{pmatrix}\right\}\text{ and}\\
\Gamma(C^1_{1,n})\;\;&=\left\{\left(\begin{pmatrix} n-1 \\ 0\end{pmatrix},\begin{pmatrix} n-1 \\ 1\end{pmatrix},\ldots,\begin{pmatrix} n-1 \\ n-1\end{pmatrix}\right)\right\}.
\end{align*}
\begin{proof} Let $n\in\N$.
First we identify $Z$ with its image $\rho_Z(Z)\sub\bigoplus_{k=1}^n H^k_n$ in $T^*(Z)$. Recall from \cite{nealrusso-Contractiveprojectionsandoperatorspaces} and \cite{NealRusso-RepresentationofcontractivelycomplementedHilbertianoperatorspacesontheFockspace} that for every $k\in\{1,\ldots,n\}$ the space $H^k_n$ is spanned by the matrices \[
b^{n,k}_i=\sum_{I\cap J=\emptyset, (I\cup J)^c=\{i\}}
\sgn(I,i,J)E_{J,I},
\]
$i=1,\ldots,n$.
Here $\sgn(I,i,J)$ is the signature of the permutation
taking\linebreak
 $(i_1,\ldots,i_{k-1},i,j_1,\ldots,j_{n-k})$ to
$(1,\ldots,n)$, when $I=\{i_1,\ldots,i_{k-1}\}$, where
$i_1<i_2<\ldots<i_{k-1}$, and $J=\{j_1,\ldots,j_{n-k}\}$ and where
$j_1<j_2<\ldots<j_{n-k}$. To compute $\Gamma(Z)$ we have do determine a grid in $\rho_Z(Z)$, which spans $\rho_Z(Z)$. From \cite{NealRusso-RepresentationofcontractivelycomplementedHilbertianoperatorspacesontheFockspace}, §1 it is known, that for every $k=1\ldots,n$ the matrices $b^{n,k}_1,\ldots,b^{n,k}_n$ are the isometric image of a rectangular grid. Thus \[\cG=\{g_i:=(b^{n,1}_i,\ldots,b_i^{n,n}):i=1,\ldots,n\}\] is a rectangular grid spanning $Z$.

One observes that the matrices $b_i^{n,k}$ can also be written as \[
b_i^{n,k}=\sum_{I\sub\{1,\ldots,n\},\atop|I|=k-1,i\not\in I}\sgn(I,i,(I\cup\{i\})^c)E_{(I\cup\{i\})^c,I},
\]
for all $i,k\in\{1,\ldots,n\}$. Therefore
\begin{align*}
b_i^{n,k}\left(b_i^{n,k}\right)^*&=\sum_{I\sub\{1,\ldots,n\},\atop|I|=k-1,i\not\in I}\sum_{J\sub\{1,\ldots,n\},\atop|J|=k-1,i\not\in J}\sgn(I,i,(I\cup\{i\})^c)\sgn(J,i,(J\cup\{i\})^c)\\
&\;\;\;\;\;\;\;\;\;\;\;\;\;\;\;\;\;\;\;\;\;\;\;\;\;\;\;\;\;\;\;\;\;\;E_{(I\cup\{i\})^c,I}E_{J,(J\cup\{i\})^c}\\
&=\sum_{I\sub\{1,\ldots,n\},\atop|I|=k-1,i\not\in I}\sgn(I,i,(I\cup\{i\})^c)^2 E_{(I\cup\{i\})^c,I}E_{I,(I\cup\{i\})^c}\\
&=\sum_{I\sub\{1,\ldots,n\},\atop|I|=k-1,i\not\in I}E_{(I\cup\{i\})^c,(I\cup\{i\})^c},
\end{align*}
which is a matrix of rank $\begin{pmatrix} n-1 \\ k-1\end{pmatrix}$, since we have that many choices for $I\sub\{1,\ldots,n\}\setminus\{i\}$, $|I|=k-1$. We get \[
g_ig_i^*=\left(\sum_{I\sub\{1,\ldots,n\},\atop|I|=0,i\not\in I}E_{(I\cup\{i\})^c,(I\cup\{i\})^c},\ldots,
\sum_{I\sub\{1,\ldots,n\},\atop|I|=n-1,i\not\in I}E_{(I\cup\{i\})^c,(I\cup\{i\})^c}\right)
\] and therefore all elements of $\Gamma(C^1_{1,n})$ lie in the same equivalence class:
\begin{align*}
[g_ig_i^*]&=\left(\begin{pmatrix} n-1 \\ 0\end{pmatrix},\begin{pmatrix} n-1 \\ 1\end{pmatrix},\ldots,\begin{pmatrix} n-1 \\ n-1\end{pmatrix}\right)\in \Z^n\\
\end{align*}
for all $i=1,\ldots,n$.
\end{proof}
\end{prop}

\subsection{Hermitian and symplectic factors}

\begin{prop}
Let $Z$ be isometric to a Cartan factor of type II with $\dim Z\geq 10$. Then  \[
\mathcal{KG}(C^2_{n})=(\Z,\N_0,\{1,\ldots,n\},\{1,\ldots,n\},\{2\}).
\]
\begin{proof}
The universal enveloping TRO of $C^2_n$ is the $C^*$-algebra $\M_n$ by \cite{BoWe1}. A grid spanning the skew-symmetric $n\times n$-matrices is \[\cG=\{g_{i,j}:=E_{i,j}-E_{j,i}:1\leq i<j\leq n\}.\] Thus $\Gamma(C^2_{n})$ is given by the equivalence classes of
\begin{align*}
g_{i,j}g_{i,j}^*&=(E_{i,j}-E_{j,i})(E_{i,j}-E_{j,i})^*\\
&=E_{i,i}+E_{j,j},
\end{align*}
for $1\leq i<j\leq n$. These are, independent of $i$ and $j$, all $\rank$ $2$ matrices.
\end{proof}
\end{prop}

\begin{prop}
If  $Z$ is $JB^*$-triple isomorphic to the finite dimensional Cartan factor $C^3_n$, then \[
\mathcal{KG}(C^3_n)=\left(\Z,\N_0,\{1,\ldots,n\},\{1,\ldots,n\},\{1,2\}\right).
\]
\begin{proof}
The universal enveloping TRO of $Z$ is by \cite{BoWe1} completely isometric to $\M_n$, thus  $K_0^{\text{\tiny JB*}}(C^3_n)=\Z$ with positive cone $\N_0$ and double-scales $\Sigma^{\text{\tiny JB*}}_\mathcal{L}(Z)=\Sigma^{\text{\tiny JB*}}_\mathcal{R}(Z)=\{1,\ldots,n\}$. The Cartan factor $C^3_n$ is spanned by the hermitian grid\[
\left\{g_{i,j}:=E_{i,j}+E_{j,i}:1\leq i< j\leq n\right\}\cup\left\{g_{i,i}:=E_{i,i}:1\leq i\leq n\right\}.
\]
This leads to \begin{align*}
g_{i,j}g_{i,j}^*&=(E_{i,j}+E_{j,i})^2\\
&=E_{i,i}+E_{j,j}, \text{ for } 1\leq i< j\leq n,
\end{align*}
and
\[g_{i,i}g_{i,i}^*=E_{i,i}, \text{ for } 1\leq i\leq n.\]
Thus \[\Gamma(C^3_n)=\{1,2\}.\]
\end{proof}
\end{prop}

\subsection{Spin factors}

To determine the K-grid invariant of the finite dimensional spin factors we need a little preparation. We already know that if $Z$ is a spin factor with $\dim Z=k+1\geq 3$, then \[T^*(Z)=\begin{cases}

  \M_{2^{n-1}}\oplus \M_{2^{n-1}}  & \text{if }k=2n-1,\\
  \M_{2^{n}} & \text{if }k=2n.

\end{cases}\]
Therefore one can easily conclude that

\[\left(K_0^{\text{\tiny JB*}}(Z),\Sigma^{\text{\tiny JB*}}_\mathcal{L}(Z),\Sigma^{\text{\tiny JB*}}_\mathcal{R}(Z)\right)=\begin{cases}

  \left(\Z^2,\{1,\ldots,2^{n-1}\}^2,\{1,\ldots,2^{n-1}\}^2\right)   \\
   \text{if }k=2n-1,\\
  \left(\Z,\{1,\ldots,2^n\},\{1,\ldots,2^n\}\right) \\ \text{if }k=2n.
\end{cases}\]

To compute $\Gamma(Z)$ we need to determine a spin grid which spans $\rho_Z(Z)\sub T^*(Z)$. Obviously we have to distinguish between the even and odd dimensional case. From \cite{HancheOlsenStoermerJordanoperatoralgebras} a spin system is known that linearly spans $\rho_Z(Z)\sub \M_{2^n}$ as a $JB^*$-algebra, in the case that $Z$ is odd dimensional, but it is unfortunately not a spin grid. It is called the standard spin system (using the abbreviation $a^n:=\underbrace{a\otimes\ldots\otimes a}_{n\text{ times}}$, $a\in \M_2$):
\begin{align*}
                 s_0^{odd}&:=id^{n},\\
                s_1^{odd}&:=\sigma_1\otimes\id^{n-1},\\
                s_2^{odd}&:=\sigma_2\otimes\id^{n-1},\\
                s_3^{odd}&:=\sigma_3\otimes\sigma_1\otimes\id^{n-2},\\
                 s_4^{odd}&:=\sigma_3\otimes\sigma_2\otimes\id^{n-2},\\
                 s_{2l+1}^{odd}&:=\sigma_3^{l}\otimes\sigma_1\otimes\id^{n-l-1},\\
                 s_{2l+2}^{odd}&:=\sigma_3^{l}\otimes\sigma_2\otimes\id^{n-l-1},  \;\;\;\;
             \end{align*}
for $1\leq l\leq n-1.$ If we drop the last idempotent $s_{2n}$ we get a spin system which generates an even dimensional spin factor embedded in $\M_{2^n}$. However we are interested in the spin system generating $\rho_Z(Z)$ inside the universal enveloping TRO $T^*(Z)=\M_{2^{n-1}}\oplus \M_{2^{n-1}}$.

\begin{lemma}\label{evenspinsystem}
 Let $Z$ be an even dimensional spin factor with $\dim Z= 2n$ then the following idempotents define a spin system generating $\rho_Z(Z)\sub T^*(Z)=\M_{2^{n-1}}\oplus\M_{2^{n-1}}$:
\begin{align*}
                 s_0^{even}&:=\left(\id^{n-1},\id^{n-1}\right),\displaybreak[0]\\
                s_1^{even}&:=\left(\sigma_1\otimes\id^{n-2},\sigma_1\otimes\id^{n-2}\right),\displaybreak[0]\\
                s_2^{even}&:=\left(\sigma_2\otimes\id^{n-2},\sigma_2\otimes\id^{n-2}\right),\displaybreak[0]\\
                s_3^{even}&:=\left(\sigma_3\otimes\sigma_1\otimes\id^{n-3},\sigma_3\otimes\sigma_1\otimes\id^{n-3}\right),\displaybreak[0]\\
                 s_4^{even}&:=\left(\sigma_3\otimes\sigma_2\otimes\id^{n-3},\sigma_3\otimes\sigma_2\otimes\id^{n-3}\right),\displaybreak[0]\\
                 s_{2l+1}^{even}&:=\left(\sigma_3^{l}\otimes\sigma_1\otimes\id^{n-l-2},\sigma_3^{l}\otimes\sigma_1\otimes\id^{n-l-2}\right),\displaybreak[0]\\
                 s_{2l+2}^{even}&:=\left(\sigma_3^{l}\otimes\sigma_2\otimes\id^{n-l-2},\sigma_3^{l}\otimes\sigma_2\otimes\id^{n-l-2}\right),\displaybreak[0]\\ s_{2n-1}^{even}&:=\left(\sigma_3^{n-1},-\sigma_3^{n-1}\right),
                 \end{align*}
                 $\text{for } 1\leq l\leq n-2.$
\begin{proof}
This is just the image of the standard spin system under the map $\varphi:\M_{2^{n-1}}\otimes\M_2\to \M_{2^{n-1}}\oplus\M_{2^{n-1}}$, $\varphi\left(A\otimes\begin{pmatrix}\lambda_1 & 0\\ 0 & \lambda_2\end{pmatrix}\right)=(\lambda_1A,\lambda_2A)$. Restricted to the TRO-span of the standard spin system in $\M_{2^n}$, which is $\M_{2^{n-1}}\otimes \mathcal{D}$, where $\mathcal{D}$ denotes the diagonal matrices in $\M_2$, this becomes a $*$-isomorphism.
\end{proof}
\end{lemma}

Next we prove a proposition which enables us to construct a spin grid out of any given spin system. With the help of this proposition we can construct spin grids which generate the even and odd dimensional spin factors embedded in their universal enveloping TROs allowing us to compute their K-grid invariants.

\begin{prop}\label{AusSpinsystemMachGrid}
Let $S=\{\id,s_1,\ldots,s_n\}$ be a spin system. If $n$ is odd (i.e.\
the corresponding spin factor is of even dimension) we can define a
spin grid $\fG=\{u_i,\widetilde{u}_i:i=1,\ldots,n\}$ by
\[u_1:=\frac{1}{2}\lr{\id-s_1},\;\;\;\widetilde{u}_1:=-\frac{1}{2}\lr{\id+s_1}\text{ and}\]
\[u_{k+1}:=\frac{1}{2}\lr{s_{2k}+is_{2k+1}},\;\;\;\widetilde{u}_{k+1}=\frac{1}{2}\lr{s_{2k}-is_{2k+1}}\;\text{ for }
k=1,\ldots,\frac{1}{2}\lr{n-1}.\] In the case that $n$ is even, a
spin grid is given by
$\fG=\{u_i,\widetilde{u}_i:i=1,\ldots,n\}\cup\{u_0\}$ with
$u_0:=s_n$.
\begin{proof}
To prove this proposition we have to verify that all elements of $\fG$
are minimal (except $u_0$ in the case of odd dimensions) tripotents
which satisfy the spin grid axioms
$(\text{SPG1}),\ldots,(\text{SPG9})$ from \cite{BoWe1} Section 3.1. As an example we prove $(\text{SPG5})$:
\begin{description}
 \item[(a)] Let $j,k\in\{1,\ldots,n-1\}$. Using the anticommutator
 relations of the spin system we get
 \begin{align*}
  \{u_{j+1},\widetilde{u}_{k+1},\widetilde{u}_{j+1}\} &= \frac{1}{2}\lr{u_{j+1}u_{k+1}\widetilde{u}_{j+1}+\widetilde{u}_{j+1}u_{k+1}u_{j+1}}\displaybreak[0]\\
    &= \frac{1}{16}(\lr{s_{2j}+is_{2j+1}}\lr{s_{2k}+is_{2k+1}}\lr{s_{2j}-is_{2j+1}}\displaybreak[0]\\
     &\:\:\;+\lr{s_{2j}-is_{2j+1}}\lr{s_{2k}+is_{2k+1}}\lr{s_{2j}+is_{2j+1}})\displaybreak[0]\\
     &= \frac{1}{8}(s_{2j}\lr{s_{2k}+is_{2k+1}} s_{2j}-is_{2j}\lr{s_{2j}-is_{2j+1}}s_{2j+1}\displaybreak[0]\\
     &\:\:\;+is_{2j+1}\lr{s_{2k}+is_{2k+1}}s_{2k}+s_{2j+1}\lr{s_{2k}+is_{2k+1}}s_{2j+1})\displaybreak[0]\\
     &=-\frac{1}{4}\lr{s_{2k}+is_{2k+1}}\displaybreak[0]\\
     &=-\frac{1}{2}u_{k+1}.
 \end{align*}
\end{description}
\begin{description}
\item[(b)] For $j\in\{1,\ldots,n-1\}$ we have
\begin{align*}
\{u_{j+1},\widetilde{u}_1,\widetilde{u}_{j+1}\} &= -\fr{16}(\lr{s_{2j}+is_{2j+1}}\lr{\id+s_{1}}\lr{s_{2j}-is_{2j+1}}\displaybreak[0] \\
   &\;\;+\lr{s_{2j}-is_{2j+1}}\lr{\id+s_{1}}\lr{s_{2j}+is_{2j+1}})\displaybreak[0]\\
   &= -\fr{8}\lr{s_{2j}^2+s_{2j}s_1s_{2j}+s_{2j+1}^2+s_{2j+1}s_1s_{2j+1}}\displaybreak[0]\\
   &= -\frac{1}{2}u_1.
\end{align*}
\item[(c)] Similarly we get
  $\{u_1,\widetilde{u}_j,\widetilde{u}_1\}=-\frac{1}{2}u_j$ for all $j\in\{2,\ldots,n\}$.
\end{description}
\end{proof}
\end{prop}

\begin{prop}
 Let $Z$ be a finite dimensional spin factor with $\dim Z=k+1.$

  If $Z$ is of even dimension, i.e.\ $k=2n-1$, $n\geq 2$, then the K-grid invariant of $Z$ is given by
\begin{align*}
K_0^{\text{\tiny JB*}}(Z)\;\;&=\Z^2,\\
K_0^{\text{\tiny JB*}}(Z)_+&=\N_0^2,\\
\Sigma^{\text{\tiny JB*}}_\cL\left(Z\right)\;\;&=\left\{1,\ldots,2^{n-1}\right\}^2,\\
\Sigma^{\text{\tiny JB*}}_\cR\left(Z\right)\;\;&=\left\{1,\ldots,2^{n-1}\right\}^2,\\
\Gamma\left(Z\right)\;\;&=\left\{\left(2^{n-2},2^{n-2}\right),\left(2^{n-1},2^{n-1}\right)\right\}.
\end{align*}

If $Z$ is of odd dimensions, i.e.\ $k=2n$, $n\geq 2$, then $\mathcal{KG}(Z)$ has the components
\begin{align*}
K_0^{\text{\tiny JB*}}\left(Z\right)\;\;&=\Z,\\
K_0^{\text{\tiny JB*}}(Z)_+&=\N_0,\\
\Sigma_\cL^{\text{\tiny JB*}}\left(Z\right)\;\;&=\{1,\ldots,2^{n}\},\\
\Sigma_\cR^{\text{\tiny JB*}}\left(Z\right)\;\;&=\{1,\ldots,2^{n}\},\\
\Gamma\left(Z\right)\;\;&=\left\{2^{n-1}\right\}.
\end{align*}
\begin{proof}
We have to prove the statements for $\Gamma\left(Z\right)$.
Let first $\dim Z$ be odd and $\mathcal{S}=\{\id,s^{odd}_1,\ldots,s^{odd}_{2n-2}\}$ be the standard spin system in $\rho_Z\left(Z\right)$ defined as above. By Proposition \ref{AusSpinsystemMachGrid} we can construct a spin grid $\cG$ from $\mathcal{S}$ linearly spanning $\rho_Z\left(Z\right)$.
We get\[
u_1^{odd}
=\begin{pmatrix}0 & 0\\ 0& 1\end{pmatrix}\otimes\id^{n-1},\;\;\;\;\;\;\;\;\;\;\;\;
\widetilde{u}_1^{odd}
=\begin{pmatrix}-1 & 0\\ 0& 0\end{pmatrix}\otimes\id^{n-1},\]
\begin{align*}
u_{l+1}^{odd}&=\sigma_3^{l-1}\otimes
\begin{pmatrix}
0 & 0 & 0 & 0\\
0 & 0 & 0 & 1\\
1 & 0 & 0 & 0\\
0 & 0 & 0 & 0
\end{pmatrix}\otimes \id^{n-l-1}\;\;\;\text{ and}\displaybreak[0]\\
\widetilde{u}_{l+1}^{odd}
&=\sigma_3^{l-1}\otimes
\begin{pmatrix}
0 & 0 & 1 & 0\\
0 & 0 & 0 & 0\\
0 & 0 & 0 & 0\\
0 & 1 & 0 & 0
\end{pmatrix}\otimes \id^{n-l-1}.
\end{align*}
This leads us to
\[
u_1^{odd}\left(u_1^{odd}\right)^*
=\begin{pmatrix}0 & 0\\ 0& 1\end{pmatrix}\otimes\id^{n-1},\;\;\;\;\;\;\;\;\;\;\;\;
\displaybreak[0]\\
\widetilde{u}_1^{odd}\left(\widetilde{u}_1^{odd}\right)^*
=\begin{pmatrix}1 & 0\\ 0& 0\end{pmatrix}\otimes\id^{n-1},\]
\begin{align*}
u_{l+1}^{odd}\left(u_{l+1}^{odd}\right)^*
&=\id^{l-1}\otimes\begin{pmatrix}
0 & 0 & 0 & 0\\
0 & 1 & 0 & 0\\
0 & 0 & 1 & 0\\
0 & 0 & 0 & 0
\end{pmatrix}\otimes \id^{n-l-1}\text{ and }\displaybreak[0]\\
\displaybreak[0]\\
\widetilde{u}_{l+1}^{odd}\left(\widetilde{u}_{l+1}^{odd}\right)^*&=
\id^{l-1}\otimes\begin{pmatrix}
1 & 0 & 0 & 0\\
0 & 0 & 0 & 0\\
0 & 0 & 0 & 0\\
0 & 0 & 0 & 1
\end{pmatrix}\otimes \id^{n-l-1}.
\end{align*}
Since $\rank \left(A\otimes B\right)=\left(\rank A\right)\left(\rank B\right)$ we can conclude from the above, for $\dim Z=2n-1$ with $n\geq 2$, that
\[
\Gamma\left(Z\right)=\left\{2^{n-1}\right\}.\]
If $\dim Z$ is even we can deduce from the above results and Lemma \ref{evenspinsystem} that
\begin{align*}
u_1^{even}&=\left(\begin{pmatrix}0 & 0\\ 0& 1\end{pmatrix}\otimes\id^{n-2},\begin{pmatrix}0 & 0\\ 0& 1\end{pmatrix}\otimes\id^{n-2}\right),\\
\displaybreak[0]\\
\widetilde{u}_1^{even}&=\left(\begin{pmatrix}-1 & 0\\ 0& 0\end{pmatrix}\otimes\id^{n-2},\begin{pmatrix}-1 & 0\\ 0& 0\end{pmatrix}\otimes\id^{n-2}\right),\\
\displaybreak[0]\\
u_{l+1}^{even}
&=\left(\sigma_3^{l-1}\otimes
\begin{pmatrix}
0 & 0 & 0 & 0\\
0 & 0 & 0 & 1\\
1 & 0 & 0 & 0\\
0 & 0 & 0 & 0
\end{pmatrix}\otimes \id^{n-l-2},\sigma_3^{l-1}\otimes
\begin{pmatrix}
0 & 0 & 0 & 0\\
0 & 0 & 0 & 1\\
1 & 0 & 0 & 0\\
0 & 0 & 0 & 0
\end{pmatrix}\otimes \id^{n-l-2}\right),\\
\displaybreak[0]\\
\widetilde{u}_{l+1}^{even}
&=\left(\sigma_3^{l-1}\otimes
\begin{pmatrix}
0 & 0 & 1 & 0\\
0 & 0 & 0 & 0\\
0 & 0 & 0 & 0\\
0 & 1 & 0 & 0
\end{pmatrix}\otimes \id^{n-l-2},\sigma_3^{l-1}\otimes
\begin{pmatrix}
0 & 0 & 1 & 0\\
0 & 0 & 0 & 0\\
0 & 0 & 0 & 0\\
0 & 1 & 0 & 0
\end{pmatrix}\otimes \id^{n-l-2}\right).
\end{align*}
The  element $u_0^{even}$ is given by \[u_0^{even}=s_{2n-1}^{even}=\left(\sigma_3^{n-1},-\sigma_3^{n-1}\right).\]
The corresponding projections are
\begin{align*}
u_1^{even}\left(u_1^{even}\right)^*&=\left(\begin{pmatrix}0 & 0\\ 0& 1\end{pmatrix}\otimes\id^{n-2},\begin{pmatrix}0 & 0\\ 0& 1\end{pmatrix}\otimes\id^{n-2}\right),\\
\displaybreak[0]\\
\widetilde{u}_1^{even}\left(\widetilde{u}_1^{even}\right)^*&=\left(\begin{pmatrix}1 & 0\\ 0& 0\end{pmatrix}\otimes\id^{n-2},\begin{pmatrix}1 & 0\\ 0& 0\end{pmatrix}\otimes\id^{n-2}\right),\\
\displaybreak[0]\\
u_{l+1}^{even}\left(u_{l+1}^{even}\right)^*&=\left(\id^{l-1}\otimes\begin{pmatrix}
0 & 0 & 0 & 0\\
0 & 1 & 0 & 0\\
0 & 0 & 1 & 0\\
0 & 0 & 0 & 0
\end{pmatrix}\otimes \id^{n-l-2},\right.\\
&\;\;\;\;\;\left.\id^{l-1}\otimes\begin{pmatrix}
0 & 0 & 0 & 0\\
0 & 1 & 0 & 0\\
0 & 0 & 1 & 0\\
0 & 0 & 0 & 0
\end{pmatrix}\otimes \id^{n-l-2}\right),\\
\displaybreak[0]\\
\widetilde{u}_{l+1}^{even}\left(\widetilde{u}_{l+1}^{even}\right)^*&=\left(
\id^{l-1}\otimes\begin{pmatrix}
1 & 0 & 0 & 0\\
0 & 0 & 0 & 0\\
0 & 0 & 0 & 0\\
0 & 0 & 0 & 1
\end{pmatrix}\otimes \id^{n-l-2},\right.\\
&\;\;\;\;\;\left.\vphantom{\id^{l-1}\otimes\begin{pmatrix}
1 & 0 & 0 & 0\\
0 & 0 & 0 & 0\\
0 & 0 & 0 & 0\\
0 & 0 & 0 & 1
\end{pmatrix}}\id^{l-1}\otimes\begin{pmatrix}
1 & 0 & 0 & 0\\
0 & 0 & 0 & 0\\
0 & 0 & 0 & 0\\
0 & 0 & 0 & 1
\end{pmatrix}\otimes \id^{n-l-2}\right)\text{ and}\\
\displaybreak[0]\\
u_0^{even}\left(u_0^{even}\right)^*&=\left(\id^{n-1},\id^{n-1}\right).
\end{align*}
For an even dimensional spin factor we can conclude that
\[\Gamma\left(Z\right)=\left\{\left(2^{n-2},2^{n-2}\right),\left(2^{n-1},2^{n-1}\right)\right\}.\]
\end{proof}
\end{prop}

Finally we are in the position to give the announced $K$-theoretic classification of the finite dimensional $JC^*$-triple systems. We first notice that, if $Z_1$ and $Z_2$ are two isomorphic finite dimensional $JC^*$-triple, the corresponding double-scaled ordered $K_0$-groups of their universal enveloping TROs are isomorphic. This isomorphism is given by $K_0^{\text{\tiny JB*}}(\vp)$, which also induces a bijection from $\Gamma(Z_1)$ to $\Gamma(Z_2)$ by Proposition \ref{gammaZistinvariante}. Thus we can conclude that if two finite dimensional $JC^*$-triple systems are isomorphic, then their K-grid invariants are isomorphic.

\begin{thm} Let $Z_1$ and $Z_2$ be finite dimensional $JC^*$-triple systems. If $\sigma:K_0^{\text{ \tiny JB*}}(Z_1)\to K_0^{\text{ \tiny JB*}}(Z_2)$ is an isomorphism with $\sigma(\mathcal{KG}(Z_1))=\mathcal{KG}(Z_2)$, then there exists a $JB^*$-isomorphism $\vp:Z_1\to Z_2$ such that $K_0^{\text{\tiny JB*}}(\vp)=\sigma$.
\begin{proof}
Since $\sigma(\mathcal{KG}(Z_1))=\mathcal{KG}(Z_2)$, we especially know that $\sigma$ is an isomorphism of the double-scaled ordered $K_0$-groups of $T^*(Z_1)$ and $T^*(Z_2)$. Using \ref{lifting of TRO homs} we can lift $\sigma$ to a complete isometry $\vp':T^*(Z_1)\to T^*(Z_2)$ with $K_0^{\text{\tiny JB*}}(\vp')=\sigma$.
The TRO $T^*(Z_1)$ is by \cite{Smith-Finitedimensionalinjectiveoperatorspaces} the finite sum of rectangular matrix algebras $T^*(Z_1)\simeq\bigoplus_{i=1}^p\M_{n_i,m_i}$, determined by the double-scaled ordered $K_0$-group of $T^*(Z_1)$.
Now we can use the information encoded in $\Gamma(Z_1)$ and the above list of the K-grid invariants to recover which summands correspond to which Cartan factor (the list allows no ambiguities). Since the K-grid invariant is additive by Proposition \ref{kgrid inv additiv} we can recover the image $\rho_{Z_1}(Z_1)\sub T^*(Z_1)$ up to (ternary) unitary equivalence (i.e an inner TRO-automorphism). The same works for $\rho_{Z_2}(Z_2)\sub T^*(Z_2)$.

Let $\cG_1$ be a grid spanning $Z_1$ and $\cG_2':=\vp'(\cG_1)\sub T^*(Z_2)$ its image under $\vp'$. Let $\phi:T^*(Z_2)\to T^*(Z_2)$ be the TRO-isomorphism mapping the linear span of $\cG_2'$ to $\rho_{Z_2}(Z_2)$ (we construct $\phi$ by using the universal property of $T^*(Z_2)$). Since $Z_2$ is finite dimensional so is $T^*(Z_2)$ and thus $\phi$ is automatically inner and unitary equivalent to the identity. If we put $$\vp:=\rho_{Z_2}^{-1}\circ\phi\circ\vp'\circ\rho_{Z_1}:Z_1\to Z_2,$$
where $\rho_{Z_2}^{-1}:\rho_{Z_2}(Z_2)\to Z_2$ is the inverse of $\rho_{Z_2}$ restricted to its image,
then $\vp$ is a $JB^*$-isomorphism with $K_0^{\text{\tiny JB*}}(\vp)=\sigma$.
\end{proof}
\end{thm}

\bibliographystyle{alpha}
\bibliography{litarbeit}
\end{document}